\documentclass{gtart_h}  

\def\ifplaintex{\expandafter\ifx\csname documentclass\endcsname\relax}

\ifplaintex 
\else
\fi

\expandafter\ifx\csname beginpicture\endcsname\relax
\expandafter\ifx\csname documentclass\endcsname\relax
\input pictex \else
\input prepictex \input pictex \input postpictex \fi\fi

\def\gt{{\mathsurround=0pt\it $\cal G\mskip-2mu$eometry \&\ 
$\cal T\!\!$opology}}        

\def\gtp{{\mathsurround=0pt\it $\cal G\mskip-2mu$eometry \&\ 
$\cal T\!\!$opology $\cal P\!$ublications}}  

\def\lognumber#1{\def\thelognumber{#1}}
\def\volumenumber#1{\def\thevolumenumber{#1}}
\def\papernumber#1{\def\thepapernumber{#1}}
\def\volumeyear#1{\def\thevolumeyear{#1}}

\def\pagenumbers#1#2{\def\startpage{#1}\def\finishpage{#2}}
\def\published#1{\def\publishdate{#1}}
\def\proposed#1{\def\theproposer{#1}}
\def\seconded#1{\def\theseconders{#1}}
\def\received#1{\def\receiveddate{#1}}
\def\revised#1{\def\reviseddate{#1}}
\def\accepted#1{\def\accepteddate{#1}}

\let\\\par\let\thelognumber\relax
\let\thevolumenumber\relax\let\thepapernumber\relax
\let\thevolumeyear\relax\let\thesamplenumber\relax\let\startpage\relax
\let\finishpage\relax\let\publishdate\relax\let\receiveddate\relax
\let\reviseddate\relax\let\accepteddate\relax\let\theasciititle\relax
\let\theasciiauthors\relax
\let\theasciiabstract\relax
\let\theasciiemail\relax\let\theshortauthors\relax\let\theshorttitle\relax

\long\def\maketitlep{   

\count0=\startpage

\gt\hfill      
\beginpicture
\setcoordinatesystem units <0.33truein, 0.33truein> point at 2.2 0.9
\setplotsymbol ({$\cal G$})
\plotsymbolspacing=9truept
\circulararc 315 degrees from 0 1 center at 0 0
\setplotsymbol ({$\cal T$})
\circulararc 315 degrees from 1 -1 center at 1 0
\endpicture
\break
{\small\ifx\thesamplenumber\relax 
Volume \else Sample
\fi\thevolumenumber\ (\thevolumeyear)
\startpage--\finishpage\nl
Published: \publishdate}
\vglue 0.5truein plus 0.4fil minus 0.1truein

{\parskip=0pt\leftskip 0pt plus 1fil\def\\{\par\smallskip}{\ifplaintex\large
\else\Large\fi\bf\thetitle}\par\medskip}   

\vglue 0pt plus 0.1fil 

{\parskip=0pt\leftskip 0pt plus 1fil\def\\{\par}{\sc\theauthors}
\par\medskip}

\vglue 0pt plus 0.1fil 

{\small\parskip=0pt\let\newline\\
{\leftskip 0pt plus 1fil\def\\{\par}{\sl\theaddress}\par}
\expandafter\ifx\theemail\relax    
\relax\else\vglue 5pt plus 0.02fil minus 2pt\def\\{\stdspace{\rm 
and}\stdspace} 
\cl{Email:\stdspace\tt\theemail}\fi
\ifx\theurl\relax                  
\relax\else\vglue 5pt plus 0.02fil minus 2pt\def\\{\stdspace{\rm 
and}\stdspace}
\cl{URL:\stdspace\tt\theurl}\fi\par}

\vglue 7pt plus 0.3fil minus 3pt

{\bf Abstract}
\vglue 5pt plus 0.1fil minus 2pt

\theabstract

\vglue 7pt plus 0.3fil minus 3pt

{\bf AMS Classification numbers}\quad Primary:\quad \theprimaryclass

Secondary:\quad \thesecondaryclass

\vglue 5pt plus 0.3fil minus 2pt

{\bf Keywords:}\quad \thekeywords

\vglue 10pt plus 0.5fil minus 5pt

{\small  Proposed: \theproposer\hfill Received: \receiveddate\nl
Seconded: \theseconders\hfill 
\ifx\reviseddate\relax                         
Accepted: \accepteddate                        
\else
Revised: \reviseddate                          
\fi}
\eject
}       

\let\maketitlepage\maketitlep
\let\maketitle\maketitlepage

\font\phead=cmsl9 scaled 950
\font=cmsl9 scaled 1050
\font\pnum=cmbx10 scaled 913
\font\lnum=cmbx10 
\font\pfoot=cmsl9 scaled 950
\font=cmsl9 scaled 1050
\ifplaintex
\headline{\vbox to 0pt{\vskip -4.5mm\line{\small\phead\ifnum
\count0=\startpage ISSN 1364-0380 (on line)
1465-3060 (printed) \hfill {\pnum\folio}\else\ifodd\count0\def\\{ }%
\ifx\theshorttitle\relax\thetitle\else\theshorttitle\fi\hfill{\pnum\folio}
\else\def\\{ and }{\pnum\folio}\hfill\ifx\theshortauthors\relax\theauthors
\else\theshortauthors\fi\fi\fi}\vss}}
\footline{\vbox to 0pt{\vglue 0mm\line{\small\pfoot\ifnum\count0=\startpage
\copyright\ \gtp\hfill\else
\gt, Volume \thevolumenumber\ (\thevolumeyear)\hfill\fi}\vss
}}
\else
\makeatletter
\def\@oddhead{{\small\lhead\ifnum\count0=\startpage ISSN 1364-0380 (on line)
1465-3060 (printed) \hfill {\lnum\number\count0}\else\ifodd\count0
\def\\{ }\ifx\theshorttitle\relax \thetitle \else\theshorttitle\fi\hfill
{\lnum\number\count0}\else\def\\{ and }{\lnum\number\count0}
\hfill\ifx\theshortauthors\relax 
\theauthors\else\theshortauthors\fi\fi\fi}}\def\@evenhead{\@oddhead}
\def\@oddfoot{\small\lfoot\ifnum\count0=\startpage\copyright\ \gtp\hfill\else
\gt, Volume \thevolumenumber\ (\thevolumeyear)\hfill\fi}
\def\@evenfoot{\@oddfoot}
\makeatother
\fi

\newwrite\gtoutfile
\long\gdef\makeheadfile{  
{\def\\{, }\def\s{ }
\immediate\openout\gtoutfile head.xxx
\immediate\write\gtoutfile{Proxy-for: \ifx\theasciiauthors\relax
\theauthors\else\theasciiauthors\fi\s<\ifx\theasciiemail\relax\theemail\else\theasciiemail\fi>}
\immediate\write\gtoutfile{\noexpand\\}
\immediate\write\gtoutfile{Authors: \ifx\theasciiauthors\relax
\theauthors\else\theasciiauthors\fi}
{\def\\{ }\immediate\write\gtoutfile{Title: \ifx\theasciititle\relax
\thetitle\else\theasciititle\fi}}
\immediate\write\gtoutfile{Subj-class: GT or SG or MG etc}
\immediate\write\gtoutfile{MSC-class: \theprimaryclass\ifx\thesecondaryclass\relax\else, \thesecondaryclass\fi}
\immediate\write\gtoutfile{Journal-ref: Geom. Topol. \thevolumenumber
(\thevolumeyear) \startpage-\finishpage}
\immediate\write\gtoutfile{Comments: Published by Geometry and Topology at}
\immediate\write\gtoutfile{\s\s http://www.maths.warwick.ac.uk/gt/GTVol\thevolumenumber/paper\thepapernumber.abs.html}
\immediate\write\gtoutfile{\noexpand\\}
\immediate\write\gtoutfile{}
\ifx\theasciiabstract\relax
\immediate\write\gtoutfile{\theabstract}\else
\immediate\write\gtoutfile{\theasciiabstract}\fi
\immediate\write\gtoutfile{}
\immediate\write\gtoutfile{\noexpand\\}
\immediate\write\gtoutfile{}
\immediate\closeout\gtoutfile}}  

\def\maketitlepage{\maketitlep\makeheadfile}
\let\maketitle\maketitlepage

\lognumber{431}
\received{13 March 2004}
\volumenumber{8}\papernumber{34}\volumeyear{2004}
\pagenumbers{1243}{1280}
\revised{2 July 2004}
\published{24 September 2004}
\accepted{21 September 2004}

\proposed{Yasha Eliashberg}

\seconded{Leonid Polterovich, Ronald Fintushel}

\usepackage{amsmath, amssymb}

\newcommand{\Z}{\ensuremath{\mathbb{Z}}} 
\newcommand{\N}{\ensuremath{\mathbb{N}}}
\newcommand{\R}{\ensuremath{\mathbb{R}}} 
\newcommand{\C}{\ensuremath{\mathbb{C}}}

\newcommand{\cM}{\ensuremath{\mathcal{M}}} 
\newcommand{\cP}{\ensuremath{\mathcal{P}}} 
\newcommand{\cU}{\ensuremath{\mathcal{U}}}

\newcommand{\Sf}{\ensuremath{\mathsf{f}}} 
 
\newtheorem{cond}{Condition}[section]
\newtheorem{defn}{Definition}[section]

\newtheorem{lem}{Lemma}[section]
\newtheorem{cor}{Corollary}[section]
\newtheorem*{main}{Main Theorem}
\newtheorem{theo}{Theorem}[section]
\newtheorem{prop}{Proposition}[section]
\newtheorem{rem}{Remark}[section]

\begin{document} 
\title[Cylindrical contact homology of subcritical manifolds]{Cylindrical 
contact homology of subcritical\\Stein-fillable contact manifolds} 
\author{Mei-Lin Yau} 
\address{Department of Mathematics, 
Michigan State University\\East Lansing, MI 48824, USA}   
\email{yau@math.msu.edu}  

\begin{abstract}  
We use contact handle decompositions and a stabilization process 
to compute the cylindrical contact homology of a subcritical 
Stein-fillable contact manifold with vanishing first Chern class,  
and show that it is completely determined by 
the homology of a subcritical Stein-filling 
of the contact manifold. 
\end{abstract}

\maketitlepage

\section{Introduction} \label{intro}

A 1--form $\alpha$ on a $(2n-1)$--dimensional oriented manifold
$M$ is called a
{\em contact} 1--form if  it satisfies the contact condition:
\begin{equation} \label{contactform}
\alpha \wedge (d\alpha )^{n-1}\neq 0 \ \ \mbox{everywhere}.
\end{equation} 
Its kernel $\xi =\{ \alpha =0\}$ is called a (co-orientable)
{\em contact structure}. $\xi$ is a codimension 1 tangent distribution
with maximal non-integrability.
The pair $(M,\xi )$ is called a {\em contact} 
manifold. Sometimes we write $(M,\alpha )$ to stress the
contact 1--form $\alpha$ instead of the contact structure
defined by $\alpha$. Note that if $\alpha$ is a contact 1--form 
then so is $f\alpha$ for any $f\in C^{\infty}(M,\R _+)$, and 
$\ker (\alpha )=\ker(f\alpha )$. 
In this paper we assume $\xi =\ker \alpha$ to be 
{\em positive}, ie, $\alpha \wedge (d\alpha )^{n-1}>0$ is a volume 
form of $M$. Two contact manifolds $(M,\xi )$ and 
$(M',\xi ')$ are {\em
contactomorphic} if there is a diffeomorphism $\phi\co M\to M'$
such that $\phi _*\xi =\xi '$. $\phi$ is called a {\em
contactomorphism}. Contact manifolds, which include many
$S^1$--bundles and hypersurfaces of symplectic manifolds, and
eventually every 3--manifold, were first introduced in \cite{L} and
\cite{M}, and has been under study for decades.

By the contact version of Darboux's theorem, all contact 1--forms are
locally  isomorphic, which implies that there is no local invariant for a
contact  structure. Moreover, it is proved by Gray in \cite{Gr} that if
two contact  structures on a closed contact manifold are homotopic as
contact  structures, then they are isotopic as contact structures.
Therefore  there are also no local invariants of the space of contact
structures on  a closed manifold. Note that the contact condition 
(\ref{contactform}) implies that $d\alpha$ restricts to a symplectic 
structure on $\xi$. The conformal class of such symplectic structures 
is independent of the choice of a defining contact 1--form for $\xi$. 
Thus we can endow $\xi$ with a $d\alpha$--compatible almost 
complex structure and the first Chern class $c_1(\xi )$ is an 
invariant of $\xi$. 

On the other hand, there are many contact structures which are homotopic
as hyperplane distributions (hence have the same $c_1(\xi )$) 
but not homotopic as contact structures (\cite{Gi}, \cite{Ho1}, 
\cite{Ho2}, \cite{U1}, \cite{U2}, etc).
This fact makes the classification of contact structures
an interesting and challenging problem. For contact 3--manifolds, many
nice partial results have been obtained
(\cite{E1}, \cite{Gi}, \cite{Ho1}, \cite{Ho2}, \cite{KM}). But much
less is known for higher dimensional cases (\cite{Ge1}, \cite{U1}, 
\cite{U2}).

{\em Contact Homology Theory} (\cite{E4}, see also \cite{U2}, 
\cite{B}, \cite{B2}), introduced by
Y Eliashberg and H Hofer in 1994 and has been expanded into a bigger 
framework {\em Symplectic Field Theory} (\cite{EGH}, \cite{BEHWZ})
provides  Floer--Gromov--Witten type of invariants to
distinguish non-isomorphic contact structures on closed manifolds: 
A contact 1--form $\alpha$ of $M$ associates a unique vector field
$R_{\alpha}$  which satisfies
\[
\alpha (R_{\alpha}) =1,\ \ \ d\alpha (R_{\alpha},\cdot ) =0.
\]
$R_{\alpha}$ is called the {\em Reeb vector field (of $\alpha$)}. 
$(M,\xi :=\ker \alpha )$ also associates a 
symplectic manifold $(\mathbb{R}\times M,d(e^t\alpha ))$, the 
{\em symplectization} of $(M,\xi )$, whose symplectic structure 
$d(e^t\alpha )$ depends (up to an $\mathbb{R}$--invariant 
diffeomorphism of $\mathbb{R}\times M$) only on $\xi$. Then 
contact homology of $(M,\xi )$ is defined by suitably counting 
in $\mathbb{R}\times M$ $(1+s)$--punctured pseudo-holomorphic 
spheres which converges exponentially to 
{\em good} periodic Reeb trajectories 
at $t=\pm \infty$ at punctures. In some favorable cases 
(see section 2) one can 
count only pseudo-holomorphic cylinders connecting good 
{\em contractible} Reeb orbits and define {\em 
cylindrical contact homology} $HC(M,\xi )$ of $(M,\xi )$. 
In this paper we consider only the $c_1(\xi )=0$ case, then  
$HC(M,\xi )$ is graded by 
the {\em reduced Conley--Zehnder index} of Reeb orbits.  
The construction of $HC(M,\xi )$
involves choices of a contact 1--form $\alpha$ and an 
$\alpha$--admissible almost complex structure.  
Yet the resulting contact homology is independent
of all these extra choices and is truly an invariant of isotopy
classes of contact structures. Though the full strength of 
contact homology is yet to be explored, 
some interesting classification
results have been obtained in the spirit
of (cylindrical) contact homology theory(\cite{C}, \cite{U1}, \cite{U2}, 
see also \cite{EGH}, \cite{B}).

Though contact homology is meant to distinguish non-isomorphic
contact structures, itself is actually an subject of interest.
One would like to know what contact homology tells about
a contact manifold. Thus it is important to compute some
concrete examples and develop computational mechanisms
of contact homology.

This paper focuses on the computation of cylindrical contact
homology of {\em subcritical Stein-fillable} contact manifolds.
 A complex  $n$--dimensional
Stein domain $(V,J)$ is called {\em subcritical} if it admits a
proper, strictly $J$--convex Morse function 
with finitely many critical points and all critical
points have Morse index $<n$. Such a function is called {\em 
subcritical}. A contact manifold is called {\em subcritical
Stein-fillable} if it is the boundary of some subcritical Stein domain
and its contact structure is the corresponding CR--structure, ie, the
field of maximal complex tangencies. Equivalently a subcritical
Stein-fillable $(M,\xi)$ can be identified with a regular level 
set of a subcritical strictly $J$--convex function on a Stein 
manifold. From now on
we will often use the shorthand ``{\em SSFC}" for 
``subcritical Stein-fillable
contact" and simply call a subcritical Stein-fillable contact
manifold a {\em SSFC} manifold, and similarly call a Stein-fillable 
contact manifold a {\em SFC} manifold. 
In this paper we obtain the following result.

\begin{main} \label{main}
Let $(M,\xi )$ be a $(2n-1)$--dimensional SSFC
manifold with $n\geq 2$,  $c_1(\xi )|_{\pi _2(M)}=0$,
and $(V,J)$ a subcritical Stein domain such
that $\partial V=M$ and $\xi$ is the maximal complex
subbundle of $TM$. Then 
\[ 
HC_i(M,\xi )\cong \underset{m\in \mathbb{N}\cup \{ 0\} }
{\oplus}H_{2(n+m-1)-i}(V). 
\]  
\end{main}

The Main Theorem results from the fact that, roughly speaking,
counting  pseudo-holomorphic cylinders
is equivalent to counting  gradient trajectories
that connect critical points of consecutive indexes of a
Morse function of a Stein filling of $(M,\xi )$. Hence
the theorem shows that the contact homology of a
SSFC manifold $(M,\xi )$ recovers in a way the homology of a Stein
domain bounded by $M$.

Here is a brief outline of this paper: After introducing cylindrical 
contact homology in Section \ref{CCH} we study 
in section \ref{Chandle} Reeb dynamics on subcritical 
contact handles, the building block of SSFC manifolds. Global dynamics 
on $M$ is discussed in Section \ref{Reebdyna}. It is shown there 
that, since $(M,\xi )$ is subcritical, one gets enough room to maneuver 
attaching handles and hence contact 1--forms to show that contact 
homology of $(M,\xi )$ is essentially generated by Reeb orbits 
contained in cocores of contact handles. 
To compute $HC(M,\xi )$ we introduce in Section \ref{stab} $(M',\xi ')$, 
the {\em stabilization} of $(M,\xi )$. 
$(M',\xi ')$ is a SSFC manifold containing $(M,\xi )$ as a
codimension 2 contact submanifold with $M'\setminus M \cong 
V\times S^1$ a trivial $S^1$--bundle over a Stein-filling $V$ 
of $M$. 
By shaping contact handles of $(M',\xi ')$
one finds that cylindrical contact homologies of $(M,\xi )$
and $(M',\xi ')$ can be represented by the same set of generators
with degrees shifted by 2. 
In Section \ref{hmlgy} we prove 
$HC_*(M,\xi )\cong HC_{*+2}(M',\xi ')$. In Section \ref{find} we 
prove that the counting of pseudo-holomorphic cylinders in 
$(M',\xi ')$ is equivalent to the counting of gradient 
trajectories in a subcritical Stein-filling $V$ of $M$ 
and hence deduce the Main 
Theorem. To this end  we first show that for 
generic $S^1$--invariant admissible almost complex 
structure 
the linearized $\overline{\partial}$--operator at 
an $S^1$--invariant solution is surjective. This is done by 
identifying it with the corresponding surjectivity problem 
in Floer Theory. Then by applying branched covering maps on 
$M'$ and the said surjectivity result to show that up to 
contact isotopies there are only $S^1$--invariant solutions 
to be counted.

\section{Cylindrical contact homology} \label{CCH}

Before introducing the cylindrical contact homology we 
would like to give a brief account on the reduced Conley--Zehnder 
index of a contractible Reeb orbit at first. 

Let $Sp(2n)=Sp(2n, \R )$ denote the group of symplectic $2n\times
2n$--matrices.  For a path $\Phi \co [0,1]\to Sp(2n)$ a {\em
Conley--Zehnder index}  (also called $\mu$--index)
$\mu (\Phi )$ is defined in terms of crossing numbers (\cite{RS}).
Here we refer readers to \cite{RS}  for a precise definition of $\mu$
and to \cite{CZ} for the original definition
and general properties of $\mu$. We point out here that
if $\Phi '$ is a path in $Sp(2n')$ and $\Phi ''$ is a path in $Sp(2n'')$
then $\mu (\Phi '\oplus \Phi '')=\mu (\Phi ')+\mu (\Phi '')$,
here $Sp(2n')\oplus Sp(2n'')$ is identified as a subgroup of
$Sp(2n'+2n'')$ in the obvious way. The following example shows that
when $n=1$, $\mu /2$ is roughly the winding number of $\Phi$.

\medskip
{\bf Example}\qua Fix $T>0$ and $A\in sp(2)=sl(2)$. 
Consider the path $\gamma \co [0,T]\to e^{tA} \in Sp(2)$.
Then
\begin{itemize}
\item $\mu (\gamma )=0$ if $A=\begin{pmatrix} 0 & b \\ a & 0
\end{pmatrix}$ for some constants $a>0$, $b>0$;
\item if $A=\begin{pmatrix} 0 & -1 \\ 1 & 0 \end{pmatrix}$
then $\mu (\gamma )=m$, $m=2n+1$ if $n\pi <T<(n+1)\pi$,
$m=2n$ if $T=n\pi$.
\end{itemize}

For computational convenience we define the {\em reduced}
Conley--Zehnder index (also called $\overline{\mu}$--index)
of a path $\Phi$ in $Sp(2n-2)$ 
to be 
\[ 
\overline{\mu}(\Phi )=\mu (\Phi )+(n-3). 
\] 
Fix a contact 1--form $\alpha$ on a $(2n-1)$--dimensional contact
manifold $(M,\xi )$.
Let $\gamma \co [o,\tau ]\to M$ be a Reeb trajectory with
$\dot{\gamma}(t)=R_{\alpha}(\gamma (t))$. Define the {\em action}
${\cal A}(\gamma )$ of $\gamma$ to be the number
\[
T ={\cal A}(\gamma ):=\int _{\gamma}\alpha
\]
The flow $(R_{\alpha})^t$ of $R_{\alpha}$ preserves $\xi$.
Thus the linearized Reeb flow
$(R_{\alpha})^t_*$ , when restricted on $\gamma$,
defines a path of symplectic maps
\[
\Upsilon (t)=(R_{\alpha})^t_*
(\gamma (0))\co \xi |_{\gamma (0)}\to \xi |_{\gamma (t)}.
\]
When $\gamma$ is periodic with period $T$,
$\Upsilon (T)$ is called the {\em linearized Poincar\'{e}
return  map} along $\gamma$. We call $\gamma$ {\em non-degenerate}
if 1 is not an eigenvalue of $\Upsilon (T)$, {\em simple} if $\gamma$
is not a nontrivial multiple cover of another Reeb orbit. 
A contact 1--form $\alpha$ is called {\em regular}
if every (contractible) Reeb orbit of $\alpha$ is non-degenerate.  
It is well-known that generic contact 1--forms are regular. 
 If we identify
$\xi _{\gamma (T)}$ with $\mathbb{R}^{2n-2}$ then $\Upsilon (T)\in
Sp(2n-2)$ is a symplectic matrix. The eigenvalues of a symplectic matrix
comes in pairs $\rho$, $\rho ^{-1}$. 

Assume $\gamma$ is a {\em contractible} periodic Reeb trajectory with
action $T$. Let $D$ be a spanning disc of
$\gamma$  and $\Psi \co \xi |_D\to \mathbb{R}^{2n-2}
\times D$ a symplectic trivialization of $\xi$ over $D$.
Then $(\gamma ,\Phi )$ defines a
path $(\Psi \circ \Upsilon \circ \Psi ^{-1})|_{\gamma}\co [0,T]
\to Sp(2n-2)$ starting from $Id$. 
The $\mu$--index of $(\gamma ,D)$ is defined to be
\[
\mu (\gamma ,\xi ,D):=
\mu (\Psi \circ \Upsilon \circ \Psi ^{-1}|_{\gamma}),
\]
and the corresponding $\overline{\mu}$--index is
\[  \overline{\mu}(\gamma ,\xi ,D)=\mu (\gamma ,\xi ,D) + (n-3) . \]
Since $D$ is contractible, $\overline{\mu}(\gamma ,\xi ,D)$
does not depend on
$\Psi$. Let $D'$ be another spanning disc of $\gamma$. Then
\begin{equation}
\overline{\mu}(\gamma ,\xi ,D)-\overline{\mu}(\gamma ,\xi ,D')=
2c_1(A)   \label{c1A}
\end{equation}
where $c_1(A):=c_1(\xi )(A)$, 
$c_1(\xi )$ is the first Chern class of $\xi$ and $A=[D\cup
D']\in H_2(M,\Z )$. In this paper we will only consider 
$c_1(\xi )=0$ case, therefore  
$\overline{\mu}(\gamma ,\xi )=\overline{\mu}(\gamma ,\xi ,D)$
is independent of the choice of a spanning disc and is denoted 
as $\overline{\gamma}$ for notational simplicity.

For a Reeb orbit $\gamma$ we denote by $\gamma ^m$ the
$m$-th multiple of $\gamma$. Recall $\Upsilon(T)$ the 
Poincar\'{e} return map of $\gamma$. Let $n(\gamma)$ denote 
the number of real negative eigenvalues of
$\Upsilon (T)$ from the interval $(-1,0)$. $n(\gamma )$ 
does not depend on the trivialization of $\xi _{\gamma (T)}$. 

\begin{defn} 
{\rm A Reeb orbit $\sigma$ is said to be {\em good} if  
\begin{equation} \label{good} 
 \sigma\neq \gamma^{2m} \quad \text{ for any $\gamma$  with  
$n(\gamma )=$odd, $m\in\N$.}  
\end{equation}  } 
\end{defn} 

For the rest of the paper we will use the notation $\cP=\cP(\alpha)$ 
to denote the set of all good contractible Reeb orbits of $\alpha$. 
Good contractible Reeb orbits with any positive multiplicity 
are included in $\cP$ as individual elements.  
Those
contractible orbits not included in ${\cal P}$ are called {\em bad}. 
The exclusion of these bad orbits is necessary in order to define
coherent orientations of moduli spaces of pseudo-holomorphic
curves  (see section 1.9 of \cite{EGH}
for more detail).

We now consider a class of almost complex structures 
on the symplectization 
\[ Symp(M,\alpha ):=(\mathbb{R}\times M,d(e^t\alpha ) )\] 
of $(M,\xi =\ker \alpha )$. 
An almost complex structure $J\co \xi\to
\xi$ on $\xi$ is called {\em $d\alpha$--compatible} if
\begin{eqnarray*}
d\alpha (x,Jx) & > & 0\ \ \mbox{for nonzero}\ \ x\in \xi , \\
d\alpha (Jx,Jy) & = & d\alpha (x,y) \ \ \ \mbox{for}\ \  x,y\in \xi .
\end{eqnarray*}
This compatibility property dose not depend on the choice
of $\alpha$. Note that $d\alpha (\cdot ,J\cdot )$ is a
Riemannian metric on $\xi$. A
$d\alpha$--compatible $J$ can be extended uniquely
to a $d(e^t\alpha )$--compatible almost complex structure on
$Symp(M,\alpha)$,  also denoted by $J$ by the abuse
of language, such that
\[
J(\frac{\partial}{\partial t})=R_{\alpha},\ \
J(R_{\alpha})=-\frac{\partial}{\partial t}.
\]
Such $J$ is called an {\em $\alpha$--admissible} almost
complex structure on $Symp(M,\alpha )$.
Observe that the Reeb vector field $R_{\alpha}$ satisfies
$\omega (R_{\alpha}, \cdot )=-d(e^t)$, hence is the
Hamiltonian vector field of of the function
$H\co \mathbb{R}\times M\to \mathbb{R}$,
$H(t,p)=e^t$.

Fix a contact quadruple $(M,\xi ,\alpha ,J)$ so that $J$ is
$\alpha$--admissible. We assume that $\alpha$ is regular. 
Fix a spanning disk $D_{\gamma}\subset M$ 
of $\gamma$ for each $\gamma \in {\cal P}={\cal P}(\alpha )$.
Given two Reeb orbits $\gamma _-$,  
$\gamma _+$ 
we denote by ${\cal M}_J(M;\gamma _-,\gamma _+)$ 
the moduli space of maps $(\tilde{u},j)$ where 
\begin{enumerate} 
\item $j$ is an almost complex structure on 
$\dot{S}^2:=S^2\setminus\{ 0,\infty\}$ (here we identify $S^2$ with 
$\mathbb{C}\cup \{ \infty \}$); 

\item $\tilde{u}=(a,u)\co (\dot{S}^2,j)\to 
(\mathbb{R}\times M,J)$ is a proper map and is $(j,J)$--holo\-mor\-phic, 
ie, $\tilde{u}$ satisfies
$d\tilde{u}\circ j= J\circ d\tilde{u}$; 

\item $\tilde{u}$ is asymptotically cylindrical over $\gamma _-$ 
at the negative end of $\mathbb{R}\times M$ at the puncture 
$0\in S^2$; and $\tilde{u}$ is asymptotically cylindrical 
over $\gamma _+$ at the positive end of $\mathbb{R}\times M$ 
at the puncture $\infty\in S^2$;  

\item $(\tilde{u},j)\thicksim (\tilde{v},j')$ if there is a 
diffeomorphism $f\co \dot{S}^2\to \dot{S}^2$ such that 
$\tilde{v}\circ f=\tilde{u}$, $f_*j=j'$, 
and $f$ fixes all punctures. 
\end{enumerate} 

 For generic choice of $J$, 
$\cM (\gamma_-,\gamma_+)={\cal M}_J(M;\gamma 
_-,\gamma _+)$, if not empty, is a smooth manifold,  
\[
\dim \cM (\gamma_-,\gamma_+):=\overline{\gamma}_+-\overline{\gamma}_-
\] 
(recall that $c_1(\xi )=0$). 
Such a $J$ is called {\em regular}. Note that
since $J$ is $\mathbb{R}$--invariant, the $\mathbb{R}$--translation
along the $\mathbb{R}$--component of $\mathbb{R}\times M$ induces a free
$\mathbb{R}$--action on ${\cal M}(\gamma _-,\gamma _+)$. 
If $\tilde{u}=(a,u)\in {\cal M}(M;\gamma _-,\gamma _+)$ 
then $u^*d\alpha \geq 0$ pointwise. We have 
\[ 
0\leq E(\tilde{u}):=\int _{\dot{S}^2}u^*d\alpha ={\cal
A}_{\alpha} (\gamma _+)-{\cal A}_{\alpha}(\gamma _-). 
\] 
$E(\tilde{u})$ is called the $d\alpha$--energy of $\tilde{u}$. 
$E(\tilde{u})=0$ iff $\gamma _-=\gamma _+$, and in this case 
the moduli space consists of a single element 
$\mathbb{R}\times \gamma _+$.

We now proceed to define the cylindrical contact homology of 
a contact manifold $(M,\xi )$. 
For a regular contact 1--form $\alpha$ defining $\xi$ we define the 
associated {\em cylindrical contact complex}  
$C(\alpha )=\underset{k\in \mathbb{Z}}{\oplus} C_k(\alpha )$ 
to be the graded vector space over $\mathbb{Q}$ generated by 
elements of ${\cal P}={\cal P}(\alpha )$, where $C_k(\alpha )$ is the 
vector space spanned by elements $\gamma \in {\cal P}$ 
with $\overline{\gamma}=k$.

Now we fix a regular $\alpha$--admissible almost complex structure 
and define the boundary map $\partial:
C_{\ast}(\alpha )\to C_{\ast -1}(\alpha )$ 
as follows. 
Let $m(\gamma)$ denote the multiplicity of $\gamma\in{\cal P}$, 
then 
\[ 
\partial \gamma :=m(\gamma)\underset{\sigma \in {\cal P}, 
\overline{\sigma}=\overline{\gamma}-1}{\sum}n_{\gamma ,\sigma }\sigma 
\] 
where $n_{\gamma ,\sigma }$ is the algebraic number of 
elements 
of ${\cal M}(\sigma ,\gamma )/\mathbb{R}$, each 
element $C\in{\cal M}(\sigma
,\gamma )/\mathbb{R}$  is weighted by $\frac{1}{m(C)}$, 
where $m(C)$ is the
multiplicity of $C$.  Then extend $\partial$ 
$\mathbb{Q}$--linearly over $C(\alpha )$. 
Note that since $\alpha$ is regular,
$\sigma$ and $\gamma$ are
non-degenerate, ${\cal M}(\sigma ,\gamma
)/\mathbb{R}$ is  compact and hence a finite set. Moreover, 
for any $\gamma \in {\cal P}$ there are only finitely many 
$\sigma$ with ${\cal A}_{\alpha}(\sigma )<{\cal A}_{\alpha}(\gamma )$. 
Thus $\partial \gamma$ is a finite sum. 
 
We have the following theorem (see \cite{U2} and Remark 1.9.2 of
\cite{EGH}).

\begin{theo} \label{d2}
Let $(\alpha ,J)$ be a regular pair.
Then $\partial \circ \partial =0$ if $C_1(\alpha )=0$.
\end{theo}

To prove $\partial \circ \partial =0$ one wants to show that 
if a 2--dimensional moduli space 
${\cal M}(\gamma  _-,\gamma  _+)$ has nonempty 
boundary, then its boundary consists of ``broken cylinders" 
$C_1\# C_2$, where $C_1\in {\cal M}
(\gamma ,\gamma _+)/\mathbb{R}$,
$C_2\in {\cal M}(\gamma _-,\gamma )/\mathbb{R}$
for some $\gamma \in {\cal P}$ with
$\overline{\gamma} =\overline{\gamma  _+}-1$. If this is not 
true then the boundary of ${\cal M}(\gamma _-,\gamma _+)$ 
will involve holomorphic curves with more than one negative ends.  
Such curves are elements of some {\em 1--dimensional} moduli 
space ${\cal M}(\gamma _-,\gamma _1,\cdots, 
\gamma _j;\gamma _+)$ with $j\geq 1$, and $\gamma _-$, 
$\gamma _1$,...,$\gamma _j$ are Reeb orbits that form the 
negative ends of the holomorphic curves. But 
\[ 
\dim {\cal M}(\gamma _-,\gamma _1,\cdots, 
\gamma _j;\gamma _+)=\overline{\gamma}_+ -\overline{\gamma}_-
-\sum _{\nu =1}^j\overline{\gamma}_{\nu}= 2-
\sum _{\nu =1}^j\overline{\gamma}_{\nu} 
\] 
which is less than 1 if $C_1(\alpha )=0$. So if $C_1(\alpha )=0$ 
then $\partial \circ\partial =0$. We will see later that 
every SSFC manifold with $\dim >3$ and $c_1(\xi )=0$ 
will have $C_*(\alpha )=0$ for all $*\leq 1$.

When $\partial \circ \partial =0$ we define the 
The {\em $j$-th cylindrical contact homology group} of the pair 
$(\alpha,J)$ to be
\[
HC_j(\alpha ,J):=
\ker (\partial |_{C_j(\alpha )})/\partial (C_{j+1}(\alpha )).
\]
The following theorem, analogous to its counterpart in Floer 
theory, asserts that $HC(\alpha ,J)$ is independent of
regular pairs $(\alpha ,J)$ satisfying $C_*(\alpha )=0$ for $*=-1,0,1$,
hence is an invariant of of $(M,\xi )$ (see \cite{U2}).

\begin{theo} \label{HC}
Let $(\alpha _0=f_0\alpha ,J_0)$,
$(\alpha _1=f_1\alpha ,J_1)$ be two regular pairs.
Assume $C_i(\alpha _0)=C_i(\alpha _1)=0$ for $i=-1,0,1$.
Then there exists a natural isomorphism
\[
\phi ^{10}\co HC_i(\alpha _0,J_0)\to HC_i(\alpha _1,J_1) .
\]
If $(\alpha _2, J_2)$ is a third regular pair then
\[
\phi ^{20}=\phi ^{21}\circ \phi ^{10},\ \ \phi ^{00}=id.
\]
\end{theo}

The proof of Theorem \ref{HC} is similar to
the proof of the corresponding theorem in Floer theory. 
Here the required chain homotopies are guaranteed by the 
existence of smooth functions $f$ on $\mathbb{R}\times M$ 
such that $d(e^tf\alpha)$ is symplectic on $\mathbb{R}\times M$ and 
$e^tf\alpha$ interpolates $e^t\alpha _0$ and $e^t\alpha _1$.   
Moreover, in a similar fashion one can show that 
$HC(M,\xi _0)\cong HC(M,\xi _1)$ for isotopic contact structures 
$\xi _0$ and $\xi _1$ on $M$. 

We remark here that though the condition $C_*(\alpha )=0$ for 
$*=1,0,-1$ looks artificial, it (or similar conditions on 
$\overline{\mu}$) may impose restrictions on the topology of 
$M$ and even the type of $\xi$. For example, when $\dim M=3$ 
and $c_1(\xi )|_{\pi _2(M)}=0$ it is proved in \cite{HWZ1.5} 
that if for some $\alpha$ $\overline{\gamma}\geq 2$ for all 
contractible Reeb orbits of $\alpha$, then $\pi _2(M)=0$ 
and $\xi$ is {\em tight}, ie, there exists no embedded 
disc $D$ in $M$ such that (i) $\partial D$ is tangent to $\xi$, 
and (ii) $D$ is transversal to $\xi$ along $\partial D$  
(see for example \cite{E1}).

\section{Contact handles} \label{Chandle}

In this section we describe some basic models of contact handles.
These basic models have been provided and discussed in detail in
\cite{W}. Since contact handles are building blocks of 
SFC manifolds
we present a similar discussion here but with a focus on the 
dynamics of Reeb orbits. 

The complex
$n$--dimensional space $\mathbb{C}^n$ together with
its standard complex structure $i$ is a Stein manifold. 
Let $(x,y,z)$ be the standard coordinates of
$\mathbb{C}^n$ with respect to the decomposition
$\mathbb{C}^n=\mathbb{R}^k\times \mathbb{R}^k \times
\mathbb{C}^{n-k}$, $(k\leq n)$. $x=(x_1, ..., x_k)$,
$y=(y_1, ..., y_k)$, $z=(z_{k+1}, ..., z_n)$, $z_l=x_l+iy_l$.

\newcommand{\st}{{\rm st}}

Fix $0\leq k \leq n$ and define
\[
f_{\st}(x,y,z)=|x|^2-\frac{1}{2}|y|^2+\frac{1}{4}|z|^2 .
\]
$f_{\st}$ is a strictly $i$--convex function on $\mathbb{C}^n$.
Note that the origin $0$ is the only critical point of $f_{\st}$,
and its Morse index is $k$.

Define
\[
Y_{\st} :=(2x,-y,\frac{z}{2})=\sum_{j=1}^{k}(2x_j\frac{\partial}{\partial 
x_j}- y_j\frac{\partial}{\partial  y_j})+\sum _{l=k+1}^n\frac{1}{2}(
x_l\frac{\partial}{\partial x_l}+y_l\frac{\partial}{\partial y_l}),
\] 
$Y_{\st}=\nabla f_{\st}$, the gradient vector field of $f_{\st}$ with
respect to the Euclidean metric. Denote by $\omega _{\st}$ the 
standard symplectic structure $\sum _{j=1}^n dx_j\wedge dy_j$ 
on $\mathbb{C}^n$. We have 
$L_{Y_{\st}}\omega _{\st}=\omega _{\st}$.
$Y_{\st}$ is a complete Liouville vector field on the symplectic manifold  
$(\mathbb{C}^n,\omega _{\st})$. 

Define $\alpha _{\st}:=\omega _{\st}(Y_{\st},\cdot )$. 
$\alpha _{\st}$ restricts to a contact 1--form on ${\cal H}$ for  
any hypersurface ${\cal H}\subset \mathbb{C}$ transversal to $Y_{\st}$.  
Note that $\alpha _{\st}=-df_{\rm st}\circ i$.

Consider a function $\Sf\co \mathbb{C}^n\to
\mathbb{R}$
\begin{equation}  \label{fc}
\Sf (x,y,z)=b|x|^2 -b'|y|^2 +\sum _{k+1}^n 
\frac{|z_l|^2}{c_l^2}; \quad b>b', c_l:\text{ positive constants.}  
\end{equation}
$0$ is the only critical point of $\Sf$. 

Define $H_c:=\{ \Sf =c\}$. 
Then $H_c\pitchfork Y_{\st}$ when $c\neq 0$. Denote the punctured  
level set $H_{c=0}-\{ 0\}$ by $H_o^{\times}$. We have 
$H_o^{\times}\pitchfork Y_{\st}$. So $\alpha _{\st}$ restricts to 
a contact 1--form on each of the level sets of $\Sf$, except at 
the point of origin.  

For $c>0$ 
$H_c$ contains two special submanifolds: 
\begin{itemize} 
\item a $(2n-k-1)$--dimensional
coisotropic ellipsoid $S^+_c:=\{ |y|=0\} \cap H_c$;  and 
\item a $(2n-2k-1)$--dimensional contact
ellipsoid $S^*_c:=\{ |y|=|x|=0\} \cap H_c$. 
\end{itemize} 

When $c<0$ there is a $(k-1)$--dimensional
isotropic sphere $S^-_c:=\{ |x|=0=|z|\} \cap H_c$ 
on $H_c$. The normal bundle of $S_c^-$ has the decomposition
\[
{\cal N}(S_c^-,H_c)={\cal CSN}(S_c^-,H_c)\oplus
T^*S_c^- \oplus\mathbb{R}R_{\Sf}
\] 
where ${\cal CSN}(S_c^-,H_c)$ 
is the {\em conformal symplectic 
normal} bundle (see \cite{W}) of $S_c^- \subset H_c$, 
$R_{\Sf}$ is the Reeb vector field of $(H_c,\alpha _{\st})$. 

Note that the vector bundle $T^*S_c^-\oplus\mathbb{R}R_{\Sf}
\cong S_c^-\times \mathbb{R}^k$ is trivial and has a natural
framing $\{ \frac{\partial}{\partial x_1}, ..., 
\frac{\partial}{\partial x_k} \}$.
The vector bundle ${\cal CSN}(S_c^-,H_c)$ is also trivial (of rank
$2(n-k)$),  and has a natural framing
${\cal F}_k$=$\{ \frac{\partial}{\partial x_{k+1}},
\frac{\partial}{\partial y_{k+1}}, ..., \frac{\partial}{\partial 
x_n}, \frac{\partial}{\partial y_n} \}$.

When $c\to 0$, $S_c^+$ (or $S_c^-$) degenerates to the point $0$. 

We denote by $\xi_c$ (resp. $\xi_o$) the corresponding 
contact structure on $(H_c,\alpha _{\st})$ (resp. $(H_o^{\times},
\alpha _{\st})$). 

\begin{prop} \label{H+H-} 
Let $\Sf$ and $\bar{\Sf}$ be two quadratic (up to an addition of 
a constant) functions of index 
$k$ with respective coefficients $(b,b',c_j)$ 
and $(\bar{b},
\bar{b}',\bar{c}_j)$ satisfying conditions in (\ref{fc}). 
Then (\ref{L_Y}) implies that for any  
level sets $H_c$ of $\Sf$ and any level set 
$\bar{H}_{\bar{c}}$ of $\bar{\Sf}$ the flow of $Y_{\st}$
will induce 
\begin{enumerate} 
\item a contact isotopy between $(H_c,S_c^+,S_c^*,\xi _c)$ 
and $(\bar{H}_{\bar{c}},\bar{S}_{\bar{c}}^+,\bar{S}_{\bar{c}}^*,
\bar{\xi}_{\bar{c}})$ when  $c>0$ and 
$\bar{c}>0$; 
\item a contact isotopy between $(H_c,S_c^-,\xi _c)$ 
and $(\bar{H}_{\bar{c}},\bar{S}_{\bar{c}}^-,\bar{\xi}_{\bar{c}})$ 
when $c<0$ and $\bar{c}<0$; 
\item a contact isotopy between $(H_c - S_c^+,\xi _c)$ 
and $(\bar{H}_{\bar{c}} - \bar{S}_{\bar{c}}^-, 
\bar{\xi}_{\bar{c}})$ when $c>0$ and $\bar{c}<0$. 
\end{enumerate} 
\end{prop} 

In particular up to contact isotopy, 
the contact structures on $(H_c^{\pm},\alpha 
_{\st})$  do not depend on the choice of the coefficients $b$, 
$b'$, and 
$c_j$; the flow of $Y_{\st}$ will produce the required contact 
isotopies that even preserve submanifolds like $S^{\pm}_c$ and 
$S^*_c$. We hence have the freedom to adjust the values of $b$, 
$b'$ and $c_j$ to get Reeb vector fields with desired 
dynamical behavior. 

For notational simplicity we will from time to time 
use the following symbols: $(H_+,S_+,S_*,\xi 
_+)$   to represent $(H_c,S_c^+,S_c^*,\xi _c)$ when 
$c>0$; and $(H_-,S_-,\xi _-)$  
to represent $(H_c,S_c^-,\xi _c)$ when $c<0$. 
$(H_+, \xi _+)$ is called (a standard model of) a {\em contact
$k$--handle}. It is {\em subcritical} if $k<n$.

We now study the Hamiltonian and 
Reeb dynamics on level sets $H_{\pm}$, $H_o$, 
of $(\mathbb{C}^n,\omega _{\st}, Y_{\st},
\Sf )$. Again $\alpha _{st}$ is used as the preferred contact 1--form. 
Let $X_{\Sf}$ denote the Hamiltonian vector field of 
$\Sf$ with respect to $\omega _{\st}$, 
\begin{equation} \label{Xf} 
X_{\Sf}=\sum_{j=1}^{k}(2bx_j\frac{\partial}{\partial y_j}+
2b'y_j\frac{\partial}{\partial  x_j})+\sum _{l=k+1}^n\frac{2}{c_l^2}(
x_l\frac{\partial}{\partial y_l}-y_l\frac{\partial}{\partial x_l}). 
\end{equation}  
The Reeb vector fields on $(H^{\times}_o, \alpha _{\st})$ and 
$(H_{\pm},\alpha _{\st})$ are 
\begin{equation} \label{R&X}
R_{\Sf}:=\frac{X_{\Sf}}{\alpha_{\st}(X_{\Sf})}  
\end{equation} 
where 
\begin{equation} \label{R&X-T} 
\alpha _{\st}(X_{\Sf})=4b|x|^2+2b'|y|^2 + 
\sum _{k+1}^n |z_l|^2/{c_l}^2=3b|x|^2+3b'|y|^2+c   
\end{equation} 
is positive away from the point of origin. 
$R_{\Sf}$ and $X_{\Sf}$ have the same integral trajectories
up to a reparametrization. 
Let $\gamma \co [0,T_h]\to \mathbb{C}^n$ be a periodic 
$X_{\Sf}$--trajectory such that $\dot{\gamma}(t)=X_{\Sf}(\gamma (t))$. 
$T_h$ is then called the {\em Hamiltonian period} of $\gamma$. 
The {\em Reeb period} of $\gamma$ can be defined similarly, and is 
actually its {\em action} $\int _{\gamma}\alpha _{\st}$. 

\begin{lem} 
{\rm(i)}\qua There is no periodic Reeb trajectory on $H^{\times}_o$ and $H_-$.

{\rm(ii)}\qua On $H_+$ all periodic Reeb trajectories are contained in $S_*$. 

{\rm(iii)}\qua If $c^2_{k+1},...,c^2_{n}$ are linearly independent over 
$\mathbb{Q}$ then the Hamiltonian period of a simple periodic 
trajectory on $H_c$ with $c>0$ is $\pi c_l^2$ for some $k<l\leq n$, 
while its action is $\pi c_l^2c$. 
\end{lem}

\begin{proof}
Let $\varphi _{t}(w)=\varphi (t,w)\co \mathbb{R}\times
\mathbb{C}^n\to \mathbb{C}^n$, $\varphi (0,w)$=$id$ be
the flow of $X_{\Sf}$, $\gamma (t)=
(x(t),y(t),z(t))$ an integral trajectory of $X_{\Sf}$ with
$\dot{\gamma}(t)=X_{\Sf}(\gamma (t))$. Then on $\gamma$
\begin{alignat}{3}
\dot{x}_j & = 2b'y_j, & \quad
\dot{y}_j & = 2bx_j,  & \quad & \text{for}\quad  1\leq j\leq k  \\
\dot{x}_l & = 2y_l/c_l^2, & \quad 
\dot{y}_l & = -2x_l/c_l^2, & \quad & \text{for}\quad  
k+1\leq l\leq n
\end{alignat}
We have \def\stutt{\vrule depth 10pt width 0pt}
\[
\left\{ \begin{array}{l}
x_{j}(t)=x_{j}(0)\cosh (2\sqrt{bb'}t) 
+ y_{j}(0)\sqrt{\frac{b'}{b}} 
\sinh (2\sqrt{bb'}t)\stutt     \\
y_{j}(t)=y_{j}(0)\cosh (2\sqrt{bb'}t) 
+ x_{j}(0)\sqrt{\frac{b}{b'}} 
\sinh (2\sqrt{bb'}t) \stutt    \\
x_{l}(t)=x_{l}(0)\cos (2t/c^2_l) 
-y_{l}(0)\sin (2t/c^2_l) \stutt   \\
y_{l}(t)=y_{l}(0)\cos (2t/c^2_l) 
+x_{l}(0)\sin (2t/c^2_l)\stutt     
\end{array}
\right. 
\]
for $1\leq j \leq k$ and $k+1 \leq l \leq n$. 
$(x_j(t),y_j(t))$, if not identically zero, is hyperbolic, while  
$|z_l(t)|$ is a constant along any $\gamma$. So $\gamma$ is contained  
in $S_*$ if it is periodic. Hence (i) and (ii) are true. 

Assume $\gamma \subset S_*$.  The Hamiltonian period of the
$z_l$--component of $\gamma$ is $\pi c_l^2$. 
Hence if $c^2_{k+1},...,c^2_{n}$ are linearly
independent over $\mathbb{Q}$ then 
for any $c>0$ there are only $n-k$ simple
periodic trajectories on $H_c$. They are $\sigma _l:=
\{ |z_l|^2=c^2_lc\} $, $l=k+1,...,n$. 
The Hamiltonian period of $\sigma _l$ is $\pi c^2_l$, which is 
independent of the value of $c$, while 
the action of $\sigma _l$ is $\pi c^2_lc$. This proves (iii). 

Note that actions of simple Reeb orbits can be made as small as
we want by choosing $c^2_l$ to be small enough.
\end{proof}

\begin{theo}     \label{th:index1}
Let $H_+$ be as above. 

{\rm(I)}\qua All periodic Reeb orbits of $H_+$ are ``good" as defined 
in (\ref{good}). 

{\rm(II)}\qua If $c^2_{k+1}$,...,$c^2_{n}$ are linearly
independent over $\mathbb{Q}$ then all Reeb orbits of
$H_c$ are non-degenerate.

{\rm(III)\qua} $C_*(H_+,\alpha _{\st})=0$ if $*<2n-k-2$ or
$*-(2n-k-2)$ is odd.

{\rm(IV)\qua} Can choose (or deform) $H_+$ for any given
$m_o>0$ such that for $*\leq m_o$
\[
\text{rk}( C_*(H_+,\alpha _{\st})) =\left\{
\begin{array}{ll} 1 & \quad if\ *=2n-k-4+2i\
for\ some\ i\in \mathbb{N}; \\ 0 & \quad otherwise.
\end{array}  \right.
\]
\end{theo}

We start with the following lemma:

\begin{lem} \label{index}
Let $\gamma$ be a contractible Reeb orbit of a contact manifold
$(M,\xi )$ with contact 1--form $\alpha$. Let $D\subset M$
be a spanning disc of $\gamma$. Then
\[
\mu (\gamma, \xi |_{D}) = \mu (\gamma, (\mathbb{C}
R_{\alpha}\oplus \xi )|_D)
\]
where $R_{\alpha}$ is the Reeb vector field of $\alpha$.
\end{lem}

\begin{proof}[Proof of Lemma~\ref{index}]
Let $\varphi \co \mathbb{C}R_{\alpha }|_D
\tilde{\rightarrow}\mathbb{C}\times D$
be the vector bundle isomorphism
$\varphi (\lambda R_{\alpha}, p)$ = $(\lambda ,p)$ for
$\lambda \in \mathbb{C}$ and $p\in D$. $\varphi $ is a symplectic
trivialization of the vector bundle
$\mathbb{C}R_{\alpha }|_D$. The action of the
linearized Reeb flow on 
$\mathbb{C}R_{\alpha }|_{\gamma}$ is a
constant path (ie, a point) in $Sp(2)$ with respect to
$\varphi$. Let $\Phi$ be any symplectic trivialization of
$\xi$ over $D$, then $\varphi \oplus \Phi$ is a symplectic
trivialization of $\mathbb{C}R_{\alpha }\oplus \xi$ over $D$.
The definition of $\varphi$ implies that
$\mu (\gamma, \xi |_{D}; \Phi )$=
$\mu (\gamma, (\mathbb{C}R_{\alpha }\oplus
\xi )|_D; \varphi \oplus \Phi )$. Since the $\mu$--index
is independent of the choice of a symplectic trivialization
over a fixed spanning disc, we conclude that
$\mu (\gamma, \xi |_{D})$=$\mu (\gamma,
(\mathbb{C}R_{\alpha}\oplus
\xi )|_D)$.
\end{proof}

\begin{proof}[Proof of Theorem~\ref{th:index1}]
Let $\Phi $:=$T\mathbb{C}^{n}\tilde{\rightarrow}
\mathbb{C}^n\times \mathbb{C}^n$
be the standard trivialization of the tangent bundle
$T\mathbb{C}^n$ of  $\mathbb{C}^n$.
When restricted on $H_+$, $\Phi$ is a
trivialization  of the stabilized contact bundle
$\mathbb{C}R\oplus\xi =T_{H_+}\mathbb{C}^n$. 
Here $R=R_{\Sf}$ is the Reeb vector field of $\alpha _{\st}$ on $H_+$. 

By Lemma \ref{index} we can use
$\Phi$ to compute the $\mu$--index of any $R_{\Sf}$--orbit $\sigma$ in
$H_+$.  Moreover, $\pi _2(H_+)=0$ if $\dim H_+>3$.
When $\dim H_+$=3, $\pi _2(H_+)$ is generated by
$S_+$. The inclusion
$H_+\hookrightarrow \mathbb{C}^2$ implies that
$(\xi _+\oplus \mathbb{C}R)|_{S_+}$ and $T_{S_+}
\mathbb{C}\cong \mathbb{C}^2\times S_+$ are isomorphic
vector bundles over $S_+$.
Since $\mathbb{C}R|_{S_+}$ is a trivial bundle over
$S_+$, so is $\xi _+|_{S_+}$, which implies that $c_1(\xi _+)=0$.
Therefore
the index $\mu (\sigma ,D_{\sigma})$ 
is independent of the choice of the spanning disk
$D_{\sigma}$ in $H_+$.

Extend the linearized Reeb flow $R^t_*$ of $R$ to $T_{H_+}
\mathbb{C}^n$ by assigning $R^t_*(Y_{\st})=Y_{\st}\circ R^t$. 
We have (with respect to $\Phi$)

\[
R^t_*|_{\sigma _l}=
e^{tD}=\left( \begin{array}{ccc}
e^{tD_{1}}  &  &       \\
      & \ddots &     \\
       &   & e^{tD_{n}}
\end{array} \right)  \  \in Sp(2n)
\] 
where
\[
D_{j} = \left( \begin{array}{cc}
0   & 2b'   \\
2b   & 0   \end{array} \right)  \ \   , \  \
D_{\ell} = \frac{2}{cc^2_{\ell}}\left( \begin{array}{cc}
0   &   -1  \\
1   &    0  \end{array}  \right)  \ \in Gl(2,\mathbb{R})
\]
for $j=1,...k$, $\ell =k+1,...n$.

By easy computation one finds that for $j=1,...,k$,

\[
e^{tD_{j}} = \left( \begin{array}{cc}
\cosh 2\sqrt{2bb'}t &
\sqrt{\frac{b'}{b}}\sinh 2\sqrt{bb'}t \\
\sqrt{\frac{b}{b'}}\sinh 2\sqrt{2bb'}t &
\cosh 2\sqrt{bb'}t
\end{array}  \right)
\]
with $\det (I-e^{tD_{j}})<0$ for all $t\neq 0$.  Moreover,
each $e^{tD_{j}}$ has two real positive eigenvalues 
$\cosh 2\sqrt{bb'}t\pm \sinh 2\sqrt{bb'}t$ for all
$t$. Therefore $e^{tD_j}$ makes no rotations to the
$(x_j,y_j)$--plane, hence has no contribution to
$\mu (\sigma _l)$.

For $\ell =k+1,...,n$, we have
\[
e^{tD_{\ell}} = \left( \begin{array}{cc}
\cos \theta   &   -\sin \theta  \\
\sin \theta   &   \cos \theta
\end{array}  \right)  \ ,\quad \mbox{where}\ \ \theta =
\frac{2t}{cc^2_{\ell}}.
\]
Note that $det(I-e^{tD_{\ell}})\geq 0$, the equality holds
if and only if $t$ is an integral multiple of
$\pi cc_{\ell}^2$. 

We have for a simple Reeb orbit $\sigma$
in $H_+$: 
\begin{enumerate}
\item $\overline{\sigma}\geq (n-k-1)+2 +(n-3)=2n-k-2$. The minimum is
always achieved by some $\sigma$.
\item For $m,m'\in \mathbb{N}$ with $m>m'$,
$\overline{\sigma^m}-\overline{\sigma^{m'}}$ is a positive even number,
and is $\geq 2(m-m')$. Here $\sigma ^m$ denotes the $m$-th
multiple of $\sigma$.
\end{enumerate}
This proves Part (I) and (III) of the theorem.

If we choose to have $c^2_{k+1}$,
$c^2_{k+2}$,  ..., $c^2_{n}$ linearly independent over
$\mathbb{Q}$, then there are exactly $n-k$ 
simple Reeb orbits
$\sigma_{k+1}$,...,$\sigma_n$ as defined before. From the
computations above it is easy to see that these $\sigma _l$ and
their multiple covers are all non-degenerate. So Part (II) is
true.

To prove Part (IV)
we consider the following perturbation of $c_l$ to compute indexes.
For any (large) integer $n_o\in \mathbb{N}$, choose
$c_{k+1}$, ..., $c_{n}$ such that $c^2_{k+1}$, ..., $c^2_n$ are
linearly independent over $\mathbb{Q}$, and
$n_oc^2_n<c^2_l$ for $l=k+1$, ..., $n-1$.  Then
\begin{equation}  \label{principal}
\begin{split}
\overline{\sigma _n^m} & =
2n-k-4+2m\ \mbox{for}\ 1\leq m\leq n_o,     \\
\overline{\sigma _n^m} & \geq
2n-k-4+2m\ \mbox{for}\ m > n_o ,    \\
\overline{\sigma _l^m} & \geq
2n-k-4+2(n_o+m) \ \mbox{for}\  l>k+1,\  m\geq 1.
\end{split}
\end{equation}
Now choose $n_o$ so that $n_o>m_o$. This completes the proof of
Part (IV).
\end{proof}

\begin{rem} 
{\rm 
We call the Reeb orbit $\sigma _n$ corresponding to $c_n$
the {\em principal (Reeb) periodic trajectory} or the {\em principal
(Reeb) orbit} of $(H_+,\alpha _{\st})$ if (\ref{principal})
is satisfied. When $m_o\to \infty$ the contact complex
$C_*(H_+,\alpha _{\st})$ is essentially generated by $\sigma _n$ 
and  its positive 
multiples. We call each positive multiple of $\sigma _n$ a
{\em  principal generator} of $C_*(H_+,\alpha _{\st})$.
We see that $C(H_+,\alpha _{\st})$ stabilizes as $m_o\to \infty$,
where $C_*(H_+,\alpha _{\st})$
is a vector space of rank 1 precisely when
$*=2n-k-4+2j$ for some $j\in \mathbb{N}$, ..., otherwise it is $0$.
}
\end{rem}

We now proceed to study the local index of a 
non-periodic Reeb trajectory $\gamma$ on $H_+$. 
Recall that our contact $k$--handle is modelled on the following 
hypersurface in $\C ^n$: 
\[ 
b|x|^2 -b'|y|^2 +\sum _{k+1}^n \frac{|z_l|^2}{c_l^2}=c>0 
\] 
Since later we will see that a SSFC manifold can be constructed 
by attaching thin subcritical contact handles to a tiny tubular 
neighborhood of attaching isotropic spheres, we are 
mainly interested in the domain $U_+\subset H_+$ 
(a tubular neighborhood of the belt sphere of $H_+$) where 
$b'|y|^2\leq C$ for some constant $C>0$, $|x|^2+|z|^2$ is small. 

Recall 
the Hamiltonian vector field $X_{\Sf}$ and Reeb vector 
field $R_{\Sf}=X_{\Sf}/\alpha_{\st}(X_{\Sf})$ from 
(\ref{Xf}), (\ref{R&X}) and (\ref{R&X-T}). 
Recall that 
the standard trivialization of the tangent bundle $T\mathbb{C}^n$
induces a symplectic trivialization (with respect to $\omega _{\st}$)
$\Phi$ of the stabilized bundle $\mathbb{C}R_{\Sf}\oplus \xi _+$ 
of $\xi _+$. 

View $\gamma$ as a 
non-periodic $X_{\Sf}$--trajectory on $U_+$ with Hamiltonian 
period $T_h$, and reparametrize $\gamma$ so that 
$\dot{\gamma}(t)=X_{\Sf}(\gamma (t))$. 
Since $H_+$ is subcritical, given any positive 
number $N_o$ we can have 
\[ 
\mu (\gamma (t), \tilde{\xi},\Phi_{\st})>N_oT_h 
\] 
by thinning $U_+$, ie, by choosing to have $c_l^2$ small 
enough. 

Now we parametrize $\gamma$ as an $R_{\Sf}$--trajectory, ie, 
$\dot{\gamma}(\tau )=R_{\Sf}(\gamma (\tau ))$. Since on $U_+$ 
\[ 
\alpha_{\st}(X_{\Sf})=3b|x|^2+3b'|y|^2+c\leq 6C+4c  
\] 
 the action 
$T$ of $\gamma$ satisfies 
\[ 
T\leq C_oT_h, \qquad C_o=6C+4c .
\] 
Denote by $\psi$ the flow of $R_{\Sf}$,
and $\varphi$ the flow of $X_{\Sf}$. We have
\[
\psi  (\tau,w)=\varphi (t(\tau,w),w)
\] 
Both flows preserve $\xi$ and have  
the following relation between their linearized flows:
\begin{equation}
d\psi (\tau ,w)=d\varphi (t (\tau ,w),w)+\frac{d\varphi}{dt}(t 
(\tau ,w),w)\otimes dt \label{eq:DR+}
\end{equation}
where $\frac{d\varphi}{dt}(t,w)=X_{\Sf}\circ \varphi (t,w)$. 
The term $\frac{d\varphi}{dt}\otimes dt=X_{\Sf}\otimes dt$ is a 
path of $2n\times 2n$ matrices of rank 1 and has $\mu$--index equal 
to 0. Hence by \cite{U1} 
we have 
\[ 
|\mu (\gamma (\tau ),\tilde{\xi},\Phi_{\st})-\mu (\gamma (t),
\tilde{\xi},\Phi_{\st})|\leq 2n, 
\] 
\[ 
\mu (\gamma (\tau ), \xi ,\Phi_{\st})=\mu (\gamma (\tau ),
\tilde{\xi},\Phi_{\st}) \geq \mu (\gamma (t), \tilde{\xi},\Phi_{\st}) 
-2n >N_oT_h-2n.
\] 
So we obtain a linear inequality relating the action $T$ of a
Reeb trajectory $\gamma$ and its $\mu$--index:   
\[  
\mu (\gamma (\tau ), \xi ,\Phi_{\st})>NT-2n, \quad N=C_o^{-1}N_o
\] 
$N$ can be
made very large by thinning the {\em subcritical} handle.
Here the fact that $H_+$ is subcritical is essential to the
largeness of $N$. We summarize the above discussion about 
non-periodic trajectories in the following lemma:

\begin{lem} \label{wander1} 
Let $U_+$ be 
a tubular neighborhood of the belt sphere of a 
$(2n-1)$--dimensional subcritical 
contact $k$--handle $H_+$. 
Let $N$ be any positive number. Then by thinning $U_+$, ie, 
by choosing to have $c_l^2$ ($l=k+1,...,n$) small enough, we 
have 
\[  
\mu (\gamma , \xi ,\Phi_{\st})>NT-2n 
\] 
for any non-periodic Reeb trajectory of $\alpha_{\st}$ on $U_+$ 
with action $T$. 
\end{lem}

\section{Reeb dynamics on SSFC manifolds} \label{Reebdyna}

A closed, orientable $(2n-1)$--dimensional contact manifold $(M,\xi )$
is called {\em Stein-fillable} if there is a $2n$--dimensional Stein
domain $(V,J)$ such that 
$\partial V=M$ and $\xi $ is the
maximal complex subbundle of $TM$. $(V,J)$ 
is called a {\em Stein filling} of  $(M,\xi )$. Let $f$ be a 
strictly $J$--convex function on $V$ that extends smoothly 
to the boundary $\partial V=M$ as a constant function, then 
$\omega _f:=-dd^{J}f=-d(df \circ J)$ 
is a symplectic 2--form on $\overline{V}$
(the nondegeneracy of $\omega _f$ is ensured by the strict 
$J$--convexity of $f$),  
and $\xi$ is the kernel of the restriction of the 1--form 
$\alpha _f:=\omega (\nabla f, \cdot )$ on $TM$. Here 
$\nabla f$ is the gradient vector field of $f$ with respect to the 
Riemannian metric $g( \cdot ,\cdot ):=\omega _f(\cdot ,J\cdot )$, 
and hence is a Liouville vector field of $\omega _f$. Without loss 
of generality we may assume that $f$ is also a Morse function. 
We call a Stein-fillable $(M,\xi )$ {\em subcritical} if the 
corresponding $f$ is {\em subcritical}, ie, has no critical 
points of index $\geq n$. Notice that $(V, \omega _f, 
\nabla f, f)$ is actually an open domain of a {\em Weinstein} 
manifold. In the following we will study the Reeb dynamics of 
SSFC manifolds in the setup of Weinstein manifolds.

A {\em Weinstein} manifold is a quadruple $(W,\omega ,Y,f)$ where 
\begin{itemize} 
\item $(W,\omega )$ is a symplectic manifold, 
\item $Y$ is a complete smooth vector field on $W$, and $Y$ is 
a {\em Liouville} vector field of $(W,\omega )$, ie,
\[ 
L_Y\omega =\omega 
\] 
where $L_Y\omega$ denotes the Lie derivative of $\omega$ with 
respect to $Y$, 
\item $f$ is an exhausting Morse function on $W$, and $Y$ is 
gradient-like with respect to $f$, ie, $df(Y)>0$ except at 
critical points of $f$.  
\end{itemize} 

In this paper 
we are interested in Weinstein manifolds of {\em finite type} 
where the function $f$ has only finitely many critical points, 
and $Y$ has only finitely many zeros accordingly. 

A Weinstein manifold $(W,\omega ,Y,f)$ associates a 1--form 
$\alpha :=\omega (Y,\cdot )$ which is a primitive of $\omega$.  
Let $S\subset W$ be a hypersurface transversal to $Y$, then 
$\alpha$ restricts to a contact 1--form on $S$. Let $X$ be a 
nonvanishing vector field which span the line field ${\cal L}_S 
\subset TS$ on which $\omega$ degenerates. Then the Reeb vector 
field of $(S,\alpha |_{TS})$ is $R:=X/\alpha (X)$. 
If $S$ is also a level set of a function $h$, then $R=X_h/
\alpha (X_h)$, where $X_h$, satisfying 
$\omega (X_h, \cdot )=-dh$, is the Hamiltonian vector 
field associated to $h$. 

Let $S'\subset W$ be another hypersurface transversal to $Y$. 
Let $\zeta$ and $\zeta '$ be contact structures on $S$ and $S'$ 
defined by $\alpha$ respectively. If a reparametrized  flow of 
$Y$ induces a diffeomorphism $\varphi \co S\to S'$ then 
we have $\varphi _*\zeta =\zeta '$, hence $(S,\zeta )$ and 
$(S',\zeta ')$ are contactomorphic. This is because for any 
two smooth functions $h_1$, $h_2>0$  on $W$, we have 
\begin{equation}  
L_{h_1Y}h_2\alpha = h_1(dh_2(Y)+h_2)\alpha . \label{L_Y} 
\end{equation}  
Note that we have $\alpha (Y)=0$ by definition.

Now let $f$ be a subcritical 
strictly $J_W$--convex Morse function on a Stein manifold
$(W,J_W)$. If $W$ is of finite type,  then $\omega :=\omega _f$ 
is independent of the choice of $f$ up
to a diffeomorphism of $W$. Let $Y_f=\nabla f$.  
The quadruple $(W,\omega _f,Y_f,f)$ is then a Weinstein manifold. 
It is easy to see that a contact manifold $(M,\xi )$ is subcritical 
Stein-fillable (up to contact isotopy) 
iff it can be realized as a hypersurface in 
some subcritical $(W,\omega _f,Y_f,f)$ that is also transversal 
to $Y_f$, or equivalently, a regular level set of $f$. 
For any level set $Q$ of $f$, $\alpha _f:=\omega (Y_f,\cdot )$ 
restricted to a contact 1--form on $Q$ away from critical points 
of $f$. Moreover the Reeb vector field associated to $\alpha$ and  
the Hamiltonian vector field of $f$ have the same integral 
trajectories. 

We now proceed to study the Reeb dynamics on 
level sets of $(W,\omega _f,Y_f,f)$. First of all,  
a theorem of Eliashberg \cite{E3} states that   
one can manipulate critical points $f$ as freely as
in the smooth case. Thus  
$f$ can be assumed to has only one critical point of
index $0$  (we assume that $W$ is connected), and $f(p)<f(q)$ 
for $p,q\in Crit(f)$ if the Morse index of $p$ is less than the 
Morse index of $q$; $f(p)=f(q)$ if $p,q\in Crit(f)$ are of the same 
index. Also, following \cite{E2} a subcritical Stein
manifold $W$ of  dimension $2n\geq 4$ can reconstructed by 
attaching handles of index $<n$. 
Therefore, once a subcritical $J$--convex Morse 
function $f$ on $(W,J_W)$ (of finite type) is chosen, 
$(W,J_W)$ can be decomposed into a (finite)
union of handlebodies of subcritical indexes accordingly. 
A critical point of Morse index $k$ corresponds to exactly a 
handlebody of of index $k$, and $(W,J_W)$ can be constructed by
attaching back these handlebodies along isotropic spheres 
with specified framings in the order of handle indexes.  

A $2n$--dimensional handlebody of index $k$ is diffeomorphic to 
$D^k\times D^{2n-k}$ with boundary $S^{k-1}\times D^{2n-k}\cup 
D^k\times D^{2n-k-1}$. $S^{k-1}\times D^{2n-k}$ is to be 
glued, while $D^k\times S^{2n-k-1}$ is a contact $k$--handle. 
Thus a SSFC manifold can be constructed by attaching subcritical 
contact handles modelled on a tubular neighborhood $U_+$ of 
the belt sphere of $H_+=\{ b|x|^2-b'|y|^2+\sum
_{l=k+1}^n|z_l|^2/c^2_l=c\}$ (see \cite{CiEH} for more detail). 
We may assume that contact handles of 
$(M,\xi )$ of the same index are pairwise disjoint. 

Recall that each subcritical contact $k$--handle has only $n-k$ 
simple Reeb orbits. We may assume that all attaching $(k-1)$--spheres 
miss all the simple Reeb orbits in the middle of contact handles 
of lower indexes. Thus a SSFC manifold has two types of 
contractible Reeb orbits. {\em Type I} Reeb orbits are those contained
in  the middle of subcritical contact handles; Reeb orbits which are 
not of type I are called {\em Type II}. Type II orbits run through 
different handles.

\begin{lem}[See Lemma 3 of \cite{CiEH}] \label{bigT} 
Let $(M, \xi )$ be a SSFC manifold with a contact handle 
decomposition. Let $T$ be any positive number. 
Then up to a contact isotopy there is a defining contact 1--form 
of $(M, \xi )$ so that any Reeb trajectory which leaves a contact 
$k$--handle and return to possibly another contact $k$--handle 
has action $\geq T$. 
\end{lem} 

Here is a brief explanation of why Lemma \ref{bigT} is true. 
The attaching isotropic spheres of subcritical contact $k$--handles 
are of dimension less than $(\dim M -1)/2$, hence after isotopy 
we may assume that there are no Reeb chords connecting these 
spheres. So for any $T>0$, there is a neighborhood $\cU_k$ of these 
spheres such that any Reeb trajectory leaves 
$\cU_k$ at time $0$ will not meet 
$\cU_k$ again before time $T$. Now we glue the 
contact $k$--handles to the interior of $\cU_k$. 

By combining the proof of Proposition 1 in \cite{CiEH} and an 
estimate of $\bar{\mu}$--index based on Lemma \ref{lem:Cm} and  
the analysis on handles in the previous section 
we can derive the following lemma concerning the $\bar{\mu}$--index 
of contractible Type II Reeb orbits.

\begin{lem} \label{wander} 
Let $(M,\xi )$ be a SSFC manifold with a subcritical contact handle 
decomposition. Let $K$ be any positive number. Then up  to contact 
isotopy (by thinning handles) every Type II contractible Reeb orbit 
has $\bar{\mu}$--index greater than $K$.
\end{lem} 

The rest of this section is devoted to proving Lemma \ref{wander}. 

We have shown that over each subcritical contact handle
there is a linear inequality relating the action $T$ of a
Reeb trajectory $\gamma $
and its $\mu $--index. Namely $\mu (\gamma )
\geq N\cdot T - 2n$, $C$ is independent of $\gamma $, $N$ can be
made very large by shaping the {\em subcritical} handle.
In the following we will estimate the actual $\mu$--index of a 
contractible Type II orbit, which is related to local indexes,
the number of
times the orbit crosses different handles, the framings of the
symplectic normal bundles of the attaching isotropic spheres,
and the gluing process. We will prove that there is a linear relation
between the action of a contractible Type II orbit and
the number of times
it crosses different handles. This linear relation, together with
the said linear inequality and
the largeness of the action of any Type II orbit,
enable us to prove Lemma \ref{wander}.

We now make a digression here to prepare for the statement of the 
inequality that links all local estimates together and guarantees 
the largeness of $\mu$-- (and hence $\bar{\mu}$--) indexes of Reeb 
orbits of Type II. 

Let $S_{n}:=\{0,1,2,..., n-1\}$ be a set of $n$ ``letters", $n\geq 2$.
Define ${\cal W}_n$ to be the set of ``words" of finite length whose
letters are elements of $S_n$.

\begin{defn}
{\rm
Given $w=l_{1}\cdots l_{m}\in {\cal W}_{n}$, $l_i\in S_n$,
$w$ is called {\em jumpy } if
$l_{i}\neq l_{i+1}$ for all $1\leq i<m$.}
\end{defn}

\begin{defn}
{\rm
Given $w\in {\cal W}_{n}$, $w=l_{1}l_{2}l_{3}\cdots
l_{m}$, $w$ {\em contains a basin} or {\em
has a basin} if there is an $k$, $0<k<n$ and a subword $w'\subset w$,
$w'$=$l_{i}\cdots l_{j}$, $1<i\leq j<m$, such that
$l_{i-1}=l_{j+1}=k>l_{\nu }$ for $i\leq \nu \leq j$. $w'$ is called a
{\em basin} of $w$.}
\end{defn}

\begin{lem}[Word Lemma]
Any jumpy word $w\in {\cal W}_{n}$ of length $2^n$ must contain
a basin ($n\geq 2$). 
\end{lem}

\begin{cor}
Any jumpy word $w\in {\cal W}_{n}$ of length $m$ must contain
at least $[\frac{m}{2^n}]$ disjoint basins.
\end{cor}

\begin{proof}[Proof of the Word Lemma]
By mathematical induction. When $n=2$, there
are only two jumpy words of length 4: 0101 and 1010. Both contain
a basin ``0" with $k=1$. Hence the lemma is true for $n=2$.
Assume the lemma
holds for $n=s$. Let $w\in {\cal W}_{s+1}$ be
jumpy of length $2^{s+1}$.
If the largest letter appearing in $w$ is less than s, then it reduces
to the case $n=s$ and the statement holds again. If the letter $s$
appears in $W$ at least twice then we are done. If not, then
$w$=$w_{1}sw_{2}$. $w_{1}$, $w_{2}\in {\cal W}_{s}$
are jumpy. Observe
that one of them is of length $\geq 2^{s}$ and hence contains a basin
by assumption. This basin is also a basin of $w$. So the lemma is
true for $n=s+1$. By induction we conclude that the lemma is true for
all $n$ greater than 1.
\end{proof} 

Now let 
$k_o$ be the highest index of contact handles of $(M,\xi )$.  
We may assume that 
Lemma \ref{bigT} holds true for $(M,\xi )$. 

Let $H(k)$ denote th union of all contact $k$--handles of $(M,\xi )$. 
Define 
\[ 
{\cal H}(k):=H(k)\setminus 
\bigcup _{k'>k}H(k')  .
\] 
Let $\gamma $ be a simple Reeb orbit of Type II.
$\gamma$ associates a jumpy word $w(\gamma )\in {\cal W}_n$
constructed as follows. 

The codimension 1 boundaries of
${\cal H}(k)$, $k=0,...,k_o$,
cut $\gamma$ into $m-1$ connected curves with boundaries.
Let $\bar{k}>0$ be the maximal value of $k$ such that $\gamma \cap
{\cal H}(k)$ is not empty.
Fix a connected component of $\gamma \cap {\cal H}(\bar{k})$ and call
it $\gamma _1$. Following the Hamiltonian flow starting from 
$\gamma _1$ we write $\gamma$ as the
ordered union $\gamma =\cup _{j=1}^{{m-1}}\gamma _j$ of these
connected curves $\gamma _j$. Define for $j=1,..,m-1$ that 
$l_j:=k$ if $\gamma _j\subset{\cal H}(k)$, and $l_m:=l_1$.
Then define $w(\gamma )$
to be $w(\gamma )$:=$l_{1}l_{2}\cdots l_{m}$. $w(\gamma )$ is
well-defined up to a choice of $\gamma _1$.
Clearly $w(\gamma )$ is jumpy and contains
at least one basin.
By the Word Lemma $w(\gamma )$ contains at least
$[\frac{m}{2^{n}}]$ disjoint basins. Also, following Lemma \ref{bigT} 
we have the property that, if $l_i\cdots l_j$ is a basin of 
$w(\gamma )$ then the action of $\gamma _i\cup \cdots \cup \gamma _j$ 
must be greater that $T$ for any prescribed number $T>0$. 
We then have the following lemma: 

\begin{lem}    \label{lem:Cm}
For any $T>0$, one can modify a subcritical Stein manifold
$(W,\omega ,Y,f)$ such that
any simple Type II Hamiltonian orbit $\gamma$
with a word $w(\gamma )$ of length $m$
must have action ${\cal A}(\gamma )>C_{m}\cdot T$, where
$C_{m}=\max \{1, [\frac{m}{2^{n}}] \}$.
\end{lem}

Recall that each handle is  attached
along an isotropic $(k-1)$--sphere  
modelled on $S_-$ with a specified
framing ${\cal F}$ of the normal bundle of the sphere. According to 
Weinstein \cite{W} there is a neighborhood $U$ of $S$ in $W$, 
a neighborhood $U_-$ of $S_-$ in $\mathbb{C}^n$, and an isomorphism 
of isotropic setups 
\[ 
\phi \co (U,\omega ,Y,M\cap U,S)\to (U_-,\omega _{\st},Y_{\st},
H_-\cap U_-,S_-).  
\] 
This isomorphism identifies the 
chosen framing ${\cal F}$ of ${\cal N}(S,M)$ with the standard framing  
of ${\cal N}(S_-,H_-)$. 
The reparamertized flow of $Y_{\st}$ induces a contactomorphism 
\[ 
\eta \co U_-\setminus S_-\to U_+\setminus S_+\subset H_+, 
\] 
$H_+$ is a
contact $k$--handle as described in Proposition \ref{H+H-}.   
Via the contactomorphism 
$\eta$ the framing ${\cal F}$ induces 
a symplectic trivialization $\Phi _{\cal F}$ of $\xi$ over 
$H_+\setminus S_+$. 
Recall also that $\xi$ on $H_+$ has a symplectic trivialization 
$\Phi _{\st}$ induced by the standard symplectic trivialization of 
$T\mathbb{C}^n$. We can choose $c_l$ so that for any
Hamiltonian trajectory $\gamma$ in some ${\cal H}(p)$ 
with ${\cal A}(\gamma )=\tau$,
\begin{equation}
\mu (\gamma ,\tilde{\xi}_k,\Phi _{\cal F})>(N_1({\bf c})+N_2({\bf
c}))\tau-2n
\label{eq:Nc}
\end{equation}
where $N({\bf c})$ is a constant depending only on $c_l$'s and
is big enough so that $N({\bf c})\tau$ exceeds any of 
the Conley--Zehnder
indexes rising from the ambiguities caused by the different
symplectic trivializations $\Phi _{\cal F}$ and $\Phi _{\st}$. 
$N_2({\bf c})$ is also a very large constant.

Now
let $\gamma $ be a simple contractible Type II orbit
of $M$. Let $D\subset M$ be a closed spanning disc
of $\gamma$.
By perturbing the interior of $D$
we may assume the following condition on $D$.
\begin{cond}
Each connected component of the intersection
$D\cap {\cal H}(k)$ contains a part of $\gamma$ if it is
not empty.
\end{cond}

Let $w(\gamma )=l_1 \cdots l_m$ be the word associated to $\gamma$.
Write $D=\cup _{j=1}^{m-1}D_j$, where $D_j$ is the intersection
of $D$ with the $j$--th handle that $\gamma$ crosses.
On each $D_j$ we use the symplectic trivialization 
$\Phi _j =\Phi _{\st}$ on $\tilde{\xi}$ and
compute the local $\mu$--index
$\mu _j=\mu (\gamma _j,\tilde{\xi}, \Phi _j)$ of $\gamma _j$. We denote
by $\mu '(\gamma )$ the sum of these local indexes.
Unfortunately, $\mu '(\gamma )$ is not the Conley--Zehnder
index that we want because the local trivializations
of $\tilde{\xi}$ by $\Phi _j$ do not match up to
a symplectic trivialization of $\tilde{\xi}|_D$.
There are two types of factors
which cause mismatches of these local trivializations:

\begin{enumerate} 
\item The choices of a framing of the normal bundle of
the attaching isotropic spheres.

\item The gluing of a $k$--handle (using the flow of $Y_{\st}$).
\end{enumerate} 
Type 1 can be overcome by
choosing suitable ${\bf c}$'s (to produce large
$N_1({\bf c})$).
Type 2 happens each time $\gamma$ crosses from
one ${\cal H}(k)$ to another. The gluing map $\eta$ preserves 
contact structure $\xi$ but not contact forms. Let 
$\alpha ':=\eta^*\alpha =e^{-h}
\alpha$, and let $R$ denote the Reeb vector field of $\alpha$, then 
the Reeb vector field $R'$ of $\alpha '$ is 
\[ 
R'=e^h(R+X^{\xi}_h) 
\] 
where $X^{\xi}$ is the vector field tangent to $\xi$ and 
satisfying 
\[ 
d\alpha (X^{\xi}_h, \cdot )=-dh|_{\xi } 
\] 
Since the actual gluing take place in a thin collar of $\partial U_-$,  
we may assume that $h\sim const$ and $R'\sim e^hR$ on the collar. 
Then by mimicking the comparison of Hamiltonian flow and Reeb flow 
in the previous section, we conclude that each Type 2 error is 
bounded by $\pm 2n$.

Let $\tau _j$ be the action of
$\gamma _j$ and $\tau =\sum _{j=1}^{m-1}\tau _j$ be
the action of $\gamma$.  By Lemma~\ref{lem:Cm} we have
\[
\tau >C_mT
\]
which together with (\ref{eq:Nc}) shows
\[
\mu '(\gamma )>\sum _{j=1}^{m-1}N_1({\bf c})\cdot \tau _j
+ N_2({\bf c})\cdot C_mT-2nm .
\]
Then the actual Conley--Zehnder index
$\mu (\gamma )=\mu (\gamma ,\tilde{\xi},D)$
satisfies the inequality
\begin{equation}
\mu (\gamma )   >   \left( \sum _{j=1}^{m-1}N_1({\bf c})\cdot \tau _j
-(\mbox{error of Type 1)}\right) + N_2({\bf c})\cdot C_mT-4mn
\label{ineq:bigmu}
\end{equation}
The first term on the right hand side of (\ref{ineq:bigmu}) can
be made positive and very large by choosing suitable $c_l$
as discussed before. For the second term, recall $C_m=\max\{ 1,
[\frac{m}{2^n}]\}$, then $N_2({\bf c})\cdot C_mT-4mn$ can be very
large if we choose to have $N_2({\bf c})\gg 2^n$ (by choosing suitable
$c_l$) and $T\gg 4n$. Note that none of these $N_1({\bf c})$,
$N_2({\bf c})$ and $T$ depend on $m$ or on $\gamma$.
This completes the proof of Lemma \ref{wander}.

\section{Stabilization of $(M,\xi )$} \label{stab}

Let $(M,\xi )$ be a $(2n-1)$--dimensional contact manifold. Let 
$(W,\omega,Y,f)$ be an Weinstein manifold associated to 
$(M,\xi )$  as discussed in the previous section. We many assume that
$M=\{ f=c\}$  for some suitable constant $c$. 
 We also assume that $\dim M>3$ for the moment. 
This condition on dimension is to ensure that 
$\partial \circ \partial =0$ because $C_*(M,\alpha )=0$ for $*\leq 1$ 
when $\dim M >3$. Later we will show that $\partial \circ \partial =0$ 
also holds true when $\dim M=3$ despite of the fact that 
$C_1(M,\alpha )\neq 0$ when $\dim M=3$. 

By now we have seen that simple Type I Reeb orbits 
of a SSFC manifold $(M,\xi)$ are in one-one 
correspondence with critical points of $f$ on $W$. 
One might expect that the counting of 1--dimensional moduli of 
holomorphic cylinders here is equivalent to the counting of 
gradient trajectories of $f$ connecting critical points of 
consecutive Morse indexes. Let $p\in Crit(f)$ be of 
index $k<n$ $(2n=\dim W)$, and let 
\[ 
S_c(p)\cong\{ \Sf =c\}\cap\{ |y|=0\} \quad (\text{ see } (\ref{fc})) 
\]  
be the $(2n-k-1)$--dimensional coisotropic ellipsoid in 
the corresponding contact $k$--handle. We may identify $\gamma_p$ with 
$\{ |z_n|^2=c^2_nc\}\cap S_c(p)$. $S^1$ acts on $S_c(p)$ by rotating  
the $z_n$--plane, giving $S_c(p)$ an open book structure with 
binding $B:=S_c(p)\cap \{ z_n=0\}\cong S^{2n-k-3}$, pages 
diffeomorphic to a $(2n-k-2)$--dimensional disc $D^{2n-k-2}$, and 
$S_c(p)\setminus B\cong S^1\times D^{2n-k-2}$. 
Then by following the discussion in 
Section \ref{find} on $S^1$--invariant holomorphic curves, one 
can see that $S_c(p)\setminus\gamma_p$ is foliated by 2--dimensional 
discs bounding $\gamma_p$ and all such discs are images of 
some element of $\cM(\gamma_p)/\R$ before contact handles of higher 
indexes are attached. In particular $B$ is the parameter space 
of a connected component of $\cM(\gamma_p)/\R$. 

Suppose now that 
a contact $(k+1)$--handle (assuming $k+1<n$) 
corresponding to $q\in Crit(f)$ 
is attached along an isotropic $k$--sphere 
which intersects transversally with $S_c(p)\setminus \gamma_p$ at 
finitely many points. Assume that these intersection points 
are on distinct elements of $\cM(S_c(p);\gamma_p)/\R$. Intuition 
suggests these ``marked" (by the intersection points) 
elements in $\cM(S_c(p);\gamma_p)/\R$ may correspond to 
elements in $\cM(\gamma_q,\gamma_p)/\R$ of the resulting manifold. 
This is where we get the speculation that perhaps the counting 
of holomorphic cylinders is equivalent to the counting of 
gradient trajectories. 
Of course many works have to be done to verify 
(or disprove) such a naive speculation. 

On the other hand, if the above guess is true for  
$(M,\xi)$ viewed as a regular level set of a subcritical Weinstein 
manifold $(W,\omega,Y,f)$, then it is also true for 
$(M',\xi')$ which is the corresponding regular level set of 
$f+\kappa|z|^2\co W\times\C\to\R$ with $\kappa>0$. It turns out that 
$(M',\xi')$ has several nice features which allow 
an alternative approach of computing $HC_*(M,\xi)$ and 
establishing a relation between $HC_*(M,\xi)$ and $H_*(W)$ 
as one has expected. The rest of this section consists of 
more discussion on $(M',\xi')$, which serves as preparation for 
the next two sections.

Recall that from $(W,\omega,Y,f)$ we can define a new Weinstein 
manifold 
\[ 
(W',\omega',Y',f'):=(W\times\C,\omega 
+\omega_o,Y+Y_o,f+\kappa |z|^2) 
\] 
where $\omega _o=dx\wedge dy$ is the standard symplectic structure on 
$\C$; $Y_o:=\frac{1}{2}(x\partial_x+y\partial_y)$ is a 
Liouville vector field with respect to 
$\omega_o$; and $\kappa >0$ is a 
constant. Here $z=x+iy$ is the complex coordinate of $\C$.
$(W',\omega ',Y',f')$ is called a 
{\em stabilization} of $(W,\omega,Y,f)$. 

Consider on $W'$ the hypersurface  
\[ 
M':=\{ f'=c\} . 
\] 
$M'$ is a regular level set of $f'$. The 1--form 
$\alpha':=\omega '(Y',\cdot )$ restricts to a contact 
1--form on $M'$. Denote by $\xi '$ 
the associated contact structure on $M'$. It is easy to see 
that $(M,\xi )$ is a codimension 2 contact submanifold of $(M',\xi ')$. 
If $M$ is subcritical, then so is $M'$.  
Note that the rotation in $\mathbb{C}$ centered at $z=0$ induces 
an $S^1$--action on $M'$ that acts freely on $M'\setminus M$ 
and fixes $M$. Indeed, we can view $M'$ as an open book with 
binding $M$, pages diffeomorphic to $V$, and trivial monodromy 
$id\co V\to V$. Here $V:=\{ f<c\}\subset W$ 
is called a {\em subcritical Stein-filling} 
of $M$. 

We can smoothly embed $\mathbb{R}\times M'$ into $W'$ by identifying 
$\{ 0\} \times M'$ with $M'\subset W'$, 
and the vector field $\frac{\partial}{\partial t}$ with 
$Y'$. The image of 
$\R\times M'$ in $W'$ is then $W'_o:=W'\setminus {\cal L}$, 
where ${\cal L}$ is the closure of the stable submanifolds 
of the flow of $Y'$. The image of $\R\times M$ is 
$W_o:=W\setminus {\cal L}$.

\begin{lem}
$c_1(\xi;M)|_{\pi _2 (M)}=0$ if and only if 
$c_1(\xi ';M')|_{\pi _2 (M')}=0$.
\end{lem}

\begin{proof}
 Let $\iota \co M\hookrightarrow M'$ be the inclusion map and
$S$ represent an element of $\pi_2(M')$. Since $M'\cong M\times D^2
\cup _{M\times S^1}V\times S^1$, 
 $\dim V=2n\geq 4$ and 
$H_{2n-2}(W)=0$, $S$ can be pushed into $M$, ie, $S$
represents an element of $\pi _2 (M)$. On the other hand 
an element of $\pi _2(M)$ is also an element of $\pi _2(M')$. 
Since $\iota ^*c_1(\xi ')=c_1(\xi )$ we conclude that
$c_1(\xi;M) |_{\pi _2 (M)}=0\Leftrightarrow 
c_1(\xi ';M')|_{\pi _2(M')}=0$.
\end{proof}

Although $\alpha '$ may not be regular in the usual sense ($\alpha
'$ may have $S^1$--families of Reeb orbits), the above properties 
ensure that cylindrical contact homology of
$(\xi ',\alpha ')$ is still defined.

The Reeb vector field of $(M',\alpha ')$ is 
\[ 
R':=\frac{X_{f}+4\kappa iY_o}{\alpha (X_{f})+\kappa |z|^2}. 
\] 
We may assume that each critical point of $f$ is standard, then so are 
the critical points of $f'$ (note that 
$f$ and $f'$ have the same 
set of critical points with the same Morse indexes). Since there are
only  finitely many critical points, $\kappa$ can be chosen so that 
$\{ c_{k+1}^2, ..., c_n^2,\kappa ^{-1}\}$ is linearly independent 
over $\mathbb{Q}$ for any $\{ c_{k+1}^2, ..., c_n^2\}$ associated 
to some critical point of index $k$ of $f$. 

When $\kappa$ is much smaller than any of those 
$c_j$ then the principal Reeb orbits of $(M',\alpha ')$  
are exactly principal Reeb orbits of $(M,\alpha )$. On the 
other hand, if $\kappa$ is much bigger than any of those $c_l$ then 
the simple principal Reeb orbits of $(M',\alpha ')$ are 
\[ 
\gamma _p :=\{ (p,z)\in M' | p\in Crit(f)\} 
\] 
they are in one-one correspondence with critical points of $f$. 
Certainly the contact homology of $(M',\xi')$ does not depend on 
the choice of $\kappa$. In fact, we will show that, up to a degree 
shifting by 2, $HC(M,\xi)$ and $HC(M',\xi')$ are isomorphic. 
Thus we know about $HC_*(M,\xi)$ once we know about $HC_*(M',\xi')$.

\section{$HC_*(M,\xi ) =  HC_{*+2}(M',\xi ' )$} 
\label{hmlgy}

Recall that $Y_o=\frac{1}{2}(x\partial_x+y\partial_y)$. 
Let $Y\subset TW$ be a gradient-like vector  
field with respect to $f$ and let $b>0$ be a constant. The 
vector field $Y':=Y+bY_o\subset TW'$ is gradient-like with repsect 
to $f'$. Let ${Y'}^t$ denote the time $t$ map of the flow of $Y'$. 
We can embed $\R\times M'$ into $W'$ by identifying 
$\{ t\}\times M'\subset \R\times M'$ with ${Y'}^t(M')\subset W'$. 
In particular the vector field $\partial_t\subset T(\R\times M')$ 
is identified with $Y'$. Note that $\R\times M\subset W\times \{ 0\}$ 
under this identification.

We need to know how to count pseudo-holomorphic 
curves in $Symp(M)=\R\times M$ and 
$Symp(M')=\R\times M'$. To achieve our goal, we 
first choose a class of admissible almost complex structures. 
First of all observe that the group $S^1$ acts on $W\times 
\mathbb{C}$ by rotations on $\mathbb{C}$, sending $(p,z)$ to 
$(p,e^{i\theta}z)$ for $\theta \in S^1 \cong 
\mathbb{R}/(2\pi \mathbb{Z})$. It restricts to an $S^1$--action 
on $M'$ that fixes $M$, acts freely on $M'\setminus M$, and preserves 
$\alpha '$. 
Let $\Pi _{\mathbb{C}}$ denote the projection $W\times \mathbb{C}
\to \mathbb{C}$, and $\Pi$ the the projection $W\times \mathbb{C}
\to W$. One might expect to find an 
$\alpha '$--admissible almost complex structure which splits 
and preserves the subbundles $\Pi^*TW$ and $\Pi^*_{\C}T\C$ of 
$TW'$. This is however, not true in general. 

On the other hand, since 
\begin{equation}  \label{xi'Mdecomp} 
\xi '|_M=\xi \oplus \underline{\C}
\end{equation}  
we consider an $d\alpha '$--admissible almost complex structure 
$J'$ on $\xi'\subset TM'$ such that $J'$ 
preserves the decomposition (\ref{xi'Mdecomp}),  
and $J'=i$ when restricted to the second factor of (\ref{xi'Mdecomp}), 
here $i$ denote the standard complex structure on $\C$.   
It is easy to see that there are plenty of $d\alpha '$--compatible  
almost complex structures satisfying the above condition. 
Then we extend $J'$ to become an admissible 
almost complex structure on 
$Symp(M')\subset W'$ via the flow of ${Y'}^t$. In particular, 
$J'(Y')={Y'}^t_*R'$ on ${Y'}^t(M')$.  
We have 
\begin{equation}  \label{Ji}
J'=\begin{bmatrix}J&O\\O&i\end{bmatrix}+Q, \quad Q\to O \text{ as } 
|z|\to 0, \ z\in\C.   
\end{equation}

\begin{theo} \label{HCM=HCM'} 
$HC_*(M,\xi ;\alpha ,J)\cong HC_{*+2}(M',\xi ';\alpha ',J')$. 
\end{theo} 

\begin{proof} 
Fix any large positive integer $m_o$. $M$ can be constructed as a 
level set of the Weinstein manifold $(W,\omega ,Y,f)$ such that 
the following action condition is satisfied: Let $\gamma _p^m$ denote 
the principal generator corresponding to the critical point $p$ of $f$ 
and with multiplicity $m$. Then  
${\cal A}(\gamma _p^1)>{\cal A}(\gamma _q^m)$ when the index of $p$ is 
greater than the index of $q$ and $m\leq m_o$.

Let $\kappa$ be small enough then the principal generators 
of the contact complex $C_*(M,\alpha )$ 
are principal generators of $C_*(M',\alpha ')$. Let $\gamma _+$ and 
$\gamma _-$ be two such principal generators, then   
\[ {\cal M}:=\cM _J'(M;\gamma_-,\gamma_+)\subset {\cal M}':=
{\cal M}_{J'}(M';\gamma _-,\gamma _+). 
\] 
Given $u' \in {\cal M}'$, write 
$u'=(u_1,u_2)$ according to the splitting 
$W'=W\times \mathbb{C}$. 
Assume $u'\not\in {\cal M}$, then $u_2\not\equiv 0$. 
 The $u_2$--component of $u'$ associates two 
winding numbers (recall that $Symp(M)\subset W$) 
\begin{align*} 
n_-&:=\text{wind}(u'\cap \{ t\ll 0\} ,Symp(M)), \\ 
n_+&:=\text{wind}(u'\cap \{ t\gg 0\} ,Symp(M)). 
\end{align*} 
Since $u'(\C^*)$ and  $Symp(M)$ are 
pseudo-holomorphic submanifolds of complement dimensions, 
$u'(\C ^*)$ intersects with 
$Symp(M)$ positively at every point of the intersection 
$u'(\C ^*)\cap Symp(M)$. Thus we have 
\begin{equation} \label{n+n-} 
n_+-n_-=\# \Big( u'(\C ^*)\cap Symp(M)\Big)\geq 0. 
\end{equation}   
Write $u_2=u_2(w)$ where $w$ denotes the complex coordinate 
of $\C$. 
$u_2\co \C^*\to \C$ is a smooth function. 
Recall that we embed $Symp(M')=\R\times M'$ into 
$W'=W\times \C$ by identifying $\partial_t$ with 
$Y'=Y+bY_o$ for some constant $b>0$. The integral 
trajectories of $Y'$ perserves the value $|z|^{-b/2}$.  
Since $u'(w)$ approaches $Symp(M)\subset W\times \{ 0\}$ 
asymptotically as 
$w$ approaches either $0$ or $\infty$ we have 
\begin{gather} 
\text{ $u_2(w)$ is asymptotically holomorphic as 
$|w|\to 0$ or $\infty$},  \label{asympholo}\\ 
|u_2(w)|^{-b/2}\to 0 \quad \text{ as $|w|\to 0$ or $\infty$}.  
\label{value} 
\end{gather}    
$u_2$ can be continuously extended to 
$\C$ by 
defining $u_2(0):=0$. The extended function 
is still denoted by $u_2$ for simplicity. 

Now by (\ref{asympholo}) and (\ref{value})  we have that near $w=0$,  
\[ 
\qquad u_2(w)\sim w^{n_-} \quad \text{ for some }n_-\in\N, \  n_->b/2.  
\] 
Here $n_-$ is exactly the earlier defined winding number of $u'(\C^*)$ 
to $Symp(M)$ near $t=-\infty$. Similarly, near $w=\infty$ we 
have
\[ 
u_2(w)\sim w^{n_+} \quad \text{ for some }n_+\in\N, \ 
n_+<b/2 ,
\]  
where $n_+$ is the winding number  of $u'(\C^*)$ 
to $Symp(M)$ near $t=\infty$.
 Then we have $n_+<b/2<n_-$ and 
in particular $n_+-n_-<0$, which contradicts with (\ref{n+n-}). 
So $u_2\equiv 0$. 
We conclude that ${\cal M}={\cal M'}$. 
The degree 2 shift is an easy observation. This completes the proof. 
\end{proof}

\section{Finding $HC(M,\xi )$} \label{find}

We now proceed to compute the 
cylindrical contact homology of $(M',\xi ')$   
(again we assume that $\dim M>3$). 
This time we choose to have $\kappa \gg 1$ so that 
the principal Reeb orbits are in $M'\setminus M$ 
and they are in one-one 
correspondence with elements of $Crit(f)$. More precisely they are 
\[ 
\gamma _p:=\{ (p,z) \} \subset  M', \quad p\in Crit(f)  
\] 
with index $\overline{\gamma}_p =2n-\text{ind}(p)$. 

Let $J'$ be the same as in the previous section. To determine the 
boundary operator of the contact complex we first need to 
characterize all moduli of the form ${\cal M}'(\gamma _-,\gamma _+)$ 
of formal dimension equal to 1. 

\begin{lem}
Fix $m_o>0$ then there are contact 1--forms on $(M',\xi ')$ 
with nondegenerate Type I Reeb orbits such that if  
are of Type I with multiplicity $\leq m_o$ and if the 
formal dimension of ${\cal M}'_{J'}(\gamma _-,\gamma _+)$ is 1 
then $\gamma _{\pm}=\gamma ^m_{p_{\pm}}$ 
for some $p_{\pm}\in Crit(f)$ with $\text{ind}(p_-)=\text{ind}(p_+)+1$ 
and $\gamma _{\pm}$ have the same multiplicity $m$. 
\end{lem} 

\begin{proof} 
Let ${\cal M}'={\cal M}'_{J'}(\gamma _-,\gamma _+)$ 
be nonempty and its formal dimension is $\overline{\gamma}_+-
\overline{\gamma}_-=1$. $\gamma _+=\gamma ^{m_+}_{p_+}$ 
for some $p_+\in Crit(f)$ with 
$\text{ind}(p_+)=k_+$, $\gamma_-=\gamma ^{m_-}_{p_-}$, where 
$\gamma _{p_-}$, $m_-$, and $k_-$ are defined similarly,  
and $m_{\pm}\leq m_o$. We have 
\begin{equation*} 
\begin{split}  
1 & = 2n - k_+ +2(m_+-1) - (2n - k_-+2(m_--1)) \\ 
 & = k_--k_++2(m_+-m_-),   
\end{split}  
\end{equation*} 
therefore $k_-\neq k_+$. If $k_-<k_+$ then ${\cal A}(\gamma _+)<
{\cal A}(\gamma _-)$ which is impossible, 
so $k_->k_+$. 

Assume $k_->k_+ + 1$ then $m_+<m_-$. Note that $m_{\pm}=n_{\pm}$ 
are the winding numbers of $u'$ around 
$Symp(M)=R\times M$ near $t=\pm\infty$ 
respectively.  Since $n_+\geq n_-$ it cannot happen that 
$k_->k_+ + 1$. So we must have
$k_-=k_++1$ and hence 
$m_+=m_-$. 
\end{proof} 

By using $n_+-n_-=\# (u'(\C ^*)\cap Symp(M))$ and the positivity 
of the intersection  $u'(\C ^*)\cap Symp(M)$ we have  
the following simple lemma:

\begin{lem} 
Let ${\cal M}'={\cal M}'_{J'}(\gamma _-,\gamma _+)$ be as in the 
previous lemma and let $u'\in {\cal M}'$. Then the curve 
$u'(\C ^*)$ does not intersect $\R\times M$. 
\end{lem}

Since all curves that we are going to count are in the 
symplectization of $M'\setminus M$ we can use the diffeomorphism  
$M'\setminus M\cong V\times S^1$ to simplify the computation. 

Consider the diffeomorphism :
\[ 
\Phi \co V\times S^1\to M'\setminus M , \quad 
\Phi (x,\theta )=(x,\sqrt{\frac{c-f}{\kappa}}\theta ) 
\] 
Then  $\Phi ^*(\alpha ')=e^{-h}d\theta +\alpha$. Write 
$\Phi ^*(\alpha ')=e^{-h}(d\theta +e^h\alpha)$ then by 
using the fact that $e^{-h}(d\theta +e^h\alpha)$ is contact 
one sees that $d(e^h\alpha )$ is symplectic on $V$. So by 
abusing notations we redenote $e^h\alpha$ as $\alpha$ and 
denote $\lambda :=d\theta +\alpha$. 

$\lambda$ is a connection 1--form on the trivial principal 
bundle 
\[ V\times S^1\overset{\pi}{\to}V,
\] 
\[ d\lambda =\pi ^*\omega, \qquad \omega :=d\alpha \]  
and $h$ is a smooth Morse function on $V$ with 
$Crit(h)=Crit(f)$ and the same corresponding Morse indexes. 

With the above isomorphism understood we will from now on 
work with the contact manifold $(V\times S^1, e^{-h}\lambda )$. 
We denote the corresponding contact structure by $\xi '$.  
$\xi '$ is the horizontal lifting of $TV$ with respect to 
the connection 1--form $\lambda$. 
Let $X_h$ be the Hamiltonian  vector field of $h$ with respect to 
the symplectic 2--form 
$\omega$, ie, $\omega (X_h,\cdot )=-dh$. Then  
the Reeb vector field of $e^{-h}\lambda$ is 
\[ 
R'=\frac{e^h}{1+\alpha (X_h)}(\partial _{\theta}+X_h)  
\]  
with $1+\alpha (X_h)>0$ on $V\times S^1$. 

$S^1$ acts freely on $V\times S^1$ by rotation along 
$S^1$ fibers. Let $J'$ be an $S^1$--invariant 
$e^{-h}\lambda$--admissible almost complex structure. Since 
$\xi '$ is transversal to the fibers, $J_{\xi '}:=J'|_{\xi '}$ induces 
an $\omega$--compatible almost complex structure $\bar{J}$ on $V$ 
by 
\[ 
\bar{J}(\pi _*\eta ):=\pi _*J'\eta \quad \eta \in \xi '. 
\] 
Conversely an $\omega$--compatible almost complex 
$\bar{J}$ structure on $V$ 
induces an $S^1$--invariant  
$d\lambda$--compatible almost complex structure $J_{\xi '}$ on 
$\xi '$ which extends to be an $S^1$--invariant  
$e^{-h}\lambda$--admissible almost complex structure $J'$ 
on the symplectization of $V\times S^1$. 

Let $\gamma _{\pm}=\gamma _{p_{\pm}}^m$ for some $p_{\pm}\in 
Crit(h)$ and $u'\in {\cal M}_{J'}(\gamma _-,\gamma _+)$ 
be with $S^1$--invariant image. Let 
\[ 
u=(\bar{u}, \theta )\co \mathbb{R}\times S^1\to V\times S^1 
\]  
be the corresponding map into $V\times S^1$. 
Since $C:=u(\mathbb{R}\times S^1)$ is $S^1$--invariant and 
$\xi '$ is transversal to the $S^1$--fibers, $\xi '$ induces a 
nonsingular foliation on $C$ generated by $\xi '\cap TC$. 
We can reparametrize $u'$ so that 
\begin{equation}  \label{usut} 
\lambda (u_s)=0 ,\quad \lambda (u_t)=m. 
\end{equation}  
Here $(s,t)$ are coordinates for $\mathbb{R}\times S^1$ (so 
$z=s+it$ is the complex coordinate), and $u_s:=\dfrac{\partial u}
{\partial s}$, $u_t:=\dfrac{\partial u}{\partial t}$. 

Let $\pi _1\co T(V\times S^1)\to \xi '$ be the projection along the 
Reeb vector field $R'$. Since $u'$ is $J'$--holomorphic we have 
\begin{equation}  \label{pi-u} 
\pi _1u_s+J'\pi_1u_t=0 ,   
\end{equation}  
ie, 
\[ 
u_s-\frac{\lambda (u_s)}{1+\alpha (X_h)}(\partial _{\theta}+X_h)
+J'\Big( u_t-\frac{\lambda (u_t)}{1+\alpha (X_h)}
(\partial _{\theta}+X_h)\Big) =0 ,
\] 
which by (\ref{usut}) is reduced to 
\begin{equation}  \label{usj'ut} 
u_s +J'\Big( u_t-m\rho
({\partial _{\theta}+X_h})\Big) =0 ,
\end{equation}  
where $\rho =(1+\alpha (X_h))^{-1}$. 
Write $u_s=(\bar{u}_s,\theta _s)$, $u_t=(\bar{u}_t,\theta _t)$.
Apply $\pi _*$ to (\ref{usj'ut}) and we have 
\begin{equation} \label{CRinV} 
\bar{u}_s+\bar{J}(\bar{u}_t-m\rho X_h) =0 . 
\end{equation}  
Hence $\bar{u}$ is a finite-energy solution to the Cauchy-Riemann type 
equation as in Floer Theory with 
\[  
\bar{u}_s=m\bar{J}(\rho X_h), \quad \bar{u}_t=0 .
\] 
Note that the flow of the vector field 
$m\rho\bar{J}X_h$ is of Morse-Smale type for generic 
$\bar{J}$ hence by \cite{SZ} the linearization of (\ref{CRinV}) at 
an $S^1$--invariant solution $\bar{u}$ 
\begin{equation}  \label{Fu} 
F_{\bar{u}}\bar{\eta} =\bar{\nabla}_s\bar{\eta} +
\bar{J}\bar{\nabla}_t \bar{\eta} 
-m\bar{\nabla}_{\bar{\eta}}(\bar{J}(\rho X_h)) 
\end{equation}  
is onto for generic $\bar{J}$. Here $\bar{\nabla}$ is the 
Levi-Civita connection associated to the Riemannian metric 
$\bar{g}:=d\alpha \circ (Id\times \bar{J})$. 

Conversely an $S^1$--invariant solution $\bar{u}$ to 
(\ref{CRinV}) can be ``lifted" to an $S^1$--invariant $J'$--holomorphic 
map $u'=(a,u)$ into $Symp(V\times S^1)$ as follows. First we lift 
$\bar{u}$ to a map $u=(\bar{u},\theta )$ into $V\times S^1$ with 
$\theta =\theta (s,t)\co \mathbb{R}\times S^1\to S^1$ satisfying  
\[ 
\theta _s=-\alpha (\bar{u}_s), \quad 
\theta _t=m-\alpha (\bar{u}_t)=m. 
\] 
Such $\theta$ exists and is unique up to the addition of a constant 
rotation. The resulting map $u$ satisfies (\ref{pi-u}) and 
(\ref{usut}). 

Now solve for the function $a=a(s,t)\co \mathbb{R}\times S^{1}\to
\mathbb{R}$  which satisfies 
\begin{equation}        \label{eq:rho}
a_s=e^{-h}\lambda (u_t)=me^{-h}, \quad a_t=-e^{-h}\lambda (u_s)=0 .
  \end{equation} 
Since $\mathbb{R}\times S^1$ is a noncompact Riemann surface, there 
exists a complex-valued function $b\co \mathbb{R}\times S^1\to 
\mathbb{C}$ such that 
\[ 
b_s=e^{-h}\lambda (u_t)=me^{-h}, \quad b_t=-e^{-h}\lambda (u_s)=0 .  
\] 
$b$ is unique up to an addition of a holomorphic function. 
Write $b_1$ for the real part of $b$, and $b_2$ for the imaginary part 
of $b$. We have 
\[ 
(b_2)_{ss}+(b_2)_{tt}=-u_s(e^{-h}\lambda 
(u_s))-u_t(e^{-h}\lambda (u_t))=0 
\] 
ie, $b_2$ is harmonic, hence the imaginary part of a holomorphic 
function $\tilde{b}$. Define $a:=b-\tilde{b}$. 
$a$ is a real-valued 
function on $\mathbb{R}\times S^1$ and satisfies (\ref{eq:rho}). 
Moreover $a$ is unique up to the addition of a real constant. 
The resulting map $u'$ is $S^1$--invariant, 
$J'$--holomorphic with multiplicity $m$, 
unique up to the rotation by a constant angle and the addition  
of a real constant and satisfies 
\[
u(s,\cdot)\rightarrow \gamma _{\pm}\ \ \mbox{as}\ \  
s\rightarrow \pm \infty  , 
\]
\[
0 < \int _{\mathbb{R}\times S^1}u^*d(e^{-h}\lambda )=
\int _{\gamma _+}e^{-h}\lambda -\int _{\gamma _-}
e^{-h}\lambda =m(e^{-h(p_+)}-e^{-h(p_-)})< \infty. 
\]
Hence $u'\in {\cal M}'_{J'}(\gamma _-,\gamma _+)$.  

Now that we have establishes for each fixed $m\in \mathbb{N}$ 
(and $m\leq m_o$ for some large $m_o$) a one-one correspondence 
between (i) the (gradient-like) $\bar{J}X_h$--trajectories in $V$ that
connecting  critical points $p_{\pm}$ with
$\text{ind}(p_+)=\text{ind}(p_-)=1$  and (ii) the $S^1$--invariant
elements in  the moduli ${\cal M}'_{J'}(\gamma _-,\gamma _+)/\mathbb{R}$ 
with $\gamma _{\pm}=\gamma _{p_{\pm}}^m$. We proceed to show that 
the linearization of the operator 
$\overline{\partial}:=\overline{\partial}_{J'}$ 
at an $S^1$--invariant solution $u'\in {\cal M}'_{J'}
(\gamma _-,\gamma _+)$ is surjective for generic $J'$. 

Recall that 
\[ 
\overline{\partial}(u')=u'_s+J'u'_t. 
\] 
Let $\pi _1$, $\pi _2$ be the projections with respect to 
the orthogonal decomposition $T(Symp(V\times S^1))
\to \xi '\oplus E$, $E$ is the vector bundle spanned by 
$\partial _t$ and $R'$. 
Write $\overline{\partial}=\overline{\partial}_1+
\overline{\partial}_2$ where 
\begin{align} 
\overline{\partial}_1(u') & =\pi _1u'_s+J'\pi _1u'_t, \\ 
\overline{\partial}_2(u') & =\pi _2u'_s+J'\pi _2u'_t . 
\end{align} 
Let $D_1$, $D_2$ denote the linearizations of $\overline{\partial}_1$ 
and $\overline{\partial}_2$ at $u'$ respectively. 

\begin{lem} \label{D1onto} 
$\pi _1D_1\co W^{1,2}(\mathbb{R}\times S^1,{u'}^*\xi ')\to 
L^2(\mathbb{R}\times S^1, {u'}^*\xi ')$ is  surjective for generic
$S^1$--invariant $J'$. 
\end{lem} 

\begin{proof} 
Let $\eta \in W^{1,2}(\mathbb{R}\times S^1,{u'}^*\xi ')$. 
Let $\nabla$ be the Levi-Civita connection 
on $Symp(V\times S^1)$ with respect to the Riemannian 
metric $g':=e^{-t+h}d(e^{t-h}\lambda)\circ(Id\times J')$. 
Note that $g'|_{\xi '}=g|_{\xi '}$ where $g$ is the 
Riemannian metric on $V\times S^1$ induced by $\bar{g}$ and 
the connection 1--form $\lambda$ of the $S^1$--bundle 
$V\times S^1$ over $V$. 
Then  
\begin{equation}  
\begin{split} 
D_1(\eta ) &=\pi_1\nabla _s\eta -(\nabla_\eta \lambda)(u'_s)e^{-h}R'
+J'\pi_1\nabla _t\eta  \\ 
&\quad +\pi_1(\nabla _{\eta }J')\pi _1u'_t -  
mJ'\nabla _{\eta}(e^{-h}R'). \label{D1-1} 
\end{split}
\end{equation}  
Write $e^{-h}R'=\dfrac{\partial_{\theta}+X_h}{1+\alpha (X_h)}=\zeta 
+\partial_{\theta}$ 
with $\zeta \in \xi '$, 
\[ 
\zeta =\frac{-\alpha (X_h)\partial_{\theta}+X_h}{1+\alpha (X_h)}
=\frac{-1}{m}\cdot \pi _1u'_t  
\] 
and apply $\pi_1$ to (\ref{D1-1}) we get 
\[ 
\pi_1D_1(\eta ) =\pi_1\nabla _s\eta +J'\pi_1\nabla _t\eta  
-m\pi_1(\nabla _{\eta }J')\zeta - 
mJ'\pi _1(\nabla _{\eta}(\zeta +\partial_{\theta})). 
\] 
Let $\bar{\eta}:=\pi _*\eta$, $\bar{\zeta}:=\pi _*\zeta$. 
Since $\pi_1(\nabla _{\eta}(J'|_E))\zeta =0$ and 
\begin{align*}  
\pi _1\nabla _{\eta}\zeta &=\pi ^*\bar{\nabla}_{\bar{\eta}}\bar{\zeta} 
\qquad (\pi ^* \mbox{ means ``horizontal lifting"}) \\   
\nabla _{\eta}\partial_{\theta} &=\nabla _{\partial_{\theta}}\eta 
+[\eta ,\partial_{\theta}]=\nabla_{\partial_{\theta}}\eta -
L_{\partial_{\theta}}\eta =0 
\end{align*} 
and $\bar{\zeta}=\rho X_h$ with $\rho =(1+\alpha (X_h))^{-1}$ we have 
\begin{equation}  
\begin{split} 
\pi _*\pi _1D_1(\eta ) &=\bar{\nabla}_s\bar{\eta}+\bar{J}
\bar{\nabla}_t\bar{\eta}-m(\bar{\nabla}_{\bar{\eta}}\bar{J})
\bar{\zeta}-m\bar{J}\bar{\nabla}_{\bar{\eta}}\bar{\zeta} 
 \\ 
&=\bar{\nabla}_s\bar{\eta}+\bar{J}
\bar{\nabla}_t\bar{\eta}-m\bar{\nabla}_{\bar{\eta}}(\bar{J}(\rho X_h)) 
\\ 
&= F_{\bar{u}}(\bar{\eta}) \qquad (\mbox{see } (\ref{Fu})). 
\end{split} 
\end{equation}  
Since $F_{\bar{u}}\co W^{1,2}(\mathbb{R}\times S^1,\bar{u}^*TV)$ 
is surjective for generic $\bar{J}$ \cite{SZ} we 
conclude that $\pi _1D_1$ is surjective for generic $S^1$--invariant 
$J'$. 
\end{proof} 

Recall that $E$ denote the vector bundle over $Symp(V\times S^1)$ 
spanned by $\partial _t$ and the Reeb vector field $R'$. 

\begin{lem} \label{D2onto} 
$D_2\co W^{1,2}(\mathbb{R}\times S^1,{u'}^*E)\to 
L^2(\mathbb{R}\times S^1,{u'}^*E)$ is 
surjective for all $J'$. 
\end{lem} 

\begin{proof} 
Let $\eta \in W^{1,2}(\mathbb{R}\times S^1,{u'}^*E)$ 
and let $\nabla$ be the 
Levi-Civita connection defined in proof of Lemma \ref{D1onto}. 
Then 
\[ 
D_2(\eta )=\pi_2\nabla_s\eta +J'\pi_2\nabla_t\eta + 
m(\nabla_{\eta}J')(e^{-h}R')+mJ'\nabla_{\eta}(e^{-h}R'). 
\]  
Since $\nabla_{\eta}J'\in \mbox{End}(\xi ')$,  
$\nabla_{R'}R'=0$, $\nabla_{R'}e^{-h}=0$ and 
$\nabla_{\eta}(e^{-h}R')=0$ we have  
\[ 
D_2(\eta )=\pi_2(\nabla_s\eta +J'\nabla_t\eta ). 
\] 
Write $\eta=\eta_1\partial _t+\eta_2R'$ then 
\[ 
D_2((\eta_1,\eta_2))=(\nabla_s\eta_1-\nabla_t\eta_2, 
\nabla_s\eta_2+\nabla_t\eta_1)
\] 
is the standard $d$--bar operator on $W^{1,2}(\mathbb{R}\times 
S^1,{u'}^*E)$ with respect to the almost 
complex structure ${u'}^*(J'|_E)$
on the trivial complex  line bundle ${u'}^*E$, hence is surjective 
because $\mathbb{R}\times S^1$ is a noncompact Riemann surface. 
\end{proof} 

Lemma \ref{D1onto} and Lemma \ref{D2onto} together imply the 
following:

\begin{lem} \label{s1onto} 
Let $J'$ be an $S^1$--invariant $\alpha '$--admissible almost 
complex structure on $Symp(M')$. Then for generic $J'$ 
the linearized operator $D=D_1+D_2$ of $\overline{\partial}_{J'}$ 
is surjective at every 
$S^1$--invariant element of ${\cal M}_{J'}(\gamma _-, 
\gamma _+)$ provided that the multiplicity of $\gamma _{\pm}$ is  
small. 
\end{lem}

In the following we would like to show that up to a homotopy of  
contact 1--forms there are no elements of ${\cal M}_{J'}
(\gamma _-,\gamma _+)$ which are not $S^1$--invariant. Our proof 
is based on results from \cite{SZ} concerning finite energy solutions 
with small periods of Cauchy-Riemann type equations and the 
following construction. 

Recall that the stabilization $M'$ can be identified with 
the following hypersurface in $W\times \mathbb{C}$: 
\[ 
\{ f+\kappa |z|^2=c\} ,\quad \kappa\gg 1 \ \ \mbox{fixed},  
\] 
with contact structure defined by the 1--form 
$\alpha  '=\alpha +\alpha _o$, where 
$\alpha _o:=\omega _o(Y_o,\cdot )$, 
$\omega _o$ is the standard symplectic 2--form on $\mathbb{C}$. 
For each $k\in \mathbb{N}$ with $k>1$ the finite group 
\[ 
\mathbb{Z}_k:=\{ \theta \in S^1|\theta^k=1\} \subset S^1 
\] 
acts on $(M',\alpha ')$ via rotation in the $\mathbb{C}$--plane. 
The action preserves $\alpha '$ and induces a branched 
$k$--covering map 
\[  
\Phi _k \co M'\to M' .   
\] 
$\Phi _k\co M'\setminus M\to  M'\setminus M$ is a $k$--covering map, 
and $\Phi _k$ fixes $M$ pointwise. $\Phi _k$ induces a contact 
1--form $\alpha _k:=(\Phi_k)_*\alpha '=\alpha +\alpha _o/k$ on $M'$. Let 
$\xi _k$ denote the contact structure defined by $\alpha _k$. 
$\alpha _k$ and $\alpha '$ can be included into a smooth family of 
contact 1--forms of $M'$ 
so $\xi _k$ and $\xi '$ are isotopic as contact 
structures. Moreover $(\Phi_k)_*\xi '=\xi _k$. 

Let $J_k$ denote the $S^1$--invariant $\alpha_k$--admissible 
almost complex structure on $Symp(V\times S^1)$ such that 
\[ 
J_k=(\Phi _k)_*J'(\Phi _k^{-1})_* \quad \text { on } \xi _k.  
\] 
Clearly Lemma \ref{s1onto} also holds for generic $J_k$, $1<k\in 
\mathbb{N}$. 

For $v'_k \in {\cal M}_{J_k}(M';\gamma_-,\gamma_+)$ let 
$\tilde{v}_k\in {\cal M}_{J_k}(\gamma^k_-,\gamma^k_+)$ be the 
$k$--fold cover of $v'_k$. Then 
$\tilde{v}_k(s,t+\frac{1}{k})=\tilde{v}_k(s,t)$ and 
the pullback by $\Phi_k$ of $\tilde{v}_k$ is a $J'$--holomorphic map 
$u'_k\in {\cal M}_{J'}(M';\gamma_-,\gamma_+)$ which also satisfies 
\begin{equation} \label{periodk}  
u'_k(s,t+\frac{1}{k})=\vartheta \cdot u'_k(s,t),  
\end{equation}  
where $\vartheta$ is the generator of $\Z_k$ which represents  
the $2\pi/k$--rotation.  
Conversely if $u'_k\in {\cal M}_{J'}(M';\gamma_-,\gamma_+)$ satisfies 
(\ref{periodk}) then $\Phi _k(u'_k) \in {\cal
M}_{J_k}(M';\gamma^k_-,\gamma^k_+)$ is a $k$--cover of some 
$v'_k \in {\cal M}_{J_k}(M';\gamma_-,\gamma_+)$. 

\begin{lem} 
Assume that $\text{ind}(p_-)-\text{ind}(p_+)=1$ and 
$\gamma _{\pm}=\gamma^m_{p_{\pm}}$. 
Then there exists $k_o\in \mathbb{N}$ such that 
for all $k\geq k_o$ all elements of 
${\cal M}_{J_k}(M';\gamma_-,\gamma_+)$ are $S^1$--invariant.  
\end{lem} 

\begin{proof} 
Suppose not. Then there exists an infinite sequence $k_{\nu}$ of 
positive integers, $\underset{\nu \to \infty}\lim k_{\nu}=\infty$,  
such that for each $k_{\nu}$ the moduli space ${\cal M}_{J_{k_{\nu}}}
(M';\gamma _-,\gamma _+)$ has an element say 
$v'_{k_{\nu}}$ which  is not 
$S^1$--invariant. Let $u'_{k_{\nu}}:=\Phi 
_{k_{\nu}}^*\tilde{v}_{k_{\nu}}$  
where $\Phi _{k_{\nu}}$ is defined as before, 
$\tilde{v}_{k_{\nu}}$ is a ${k_{\nu}}$--cover of 
$v'_{k_{\nu}}$. Then $u'_{k_{\nu}}\in 
{\cal M}_{J'}(M';\gamma_-,\gamma_+)$. Since all 
$u'_{k_{\nu}}$ have the 
same contact energy 
${\cal A}_{\alpha '}(\gamma _+)-{\cal A}_{\alpha '}(\gamma
_-)$,  there is an infinite subsequence of 
$u'_{k_{\nu}}$, also denoted by $u'_{k_{\nu}}$, such that up 
to translations in $\mathbb{R}$--direction, 
$u'_{k_{\nu}}$ converge 
to a $J'$--holomorphic 
curve $u'$ as $\nu \to \infty$. $u'$ is $S^1$--invariant.  

All $u'_{k_{\nu}}$ and $u'$ have the same winding numbers 
$n_{\pm}$ around $Symp(M)$ near $t=\pm\infty$. Moreover 
we have $n_+=n_-$ because 
none of the $u'_{k_{\nu}}$'s intersect with $Symp(M)$, 
so neither does $u'$.  
$u'$ is therefore a finite union of $S^1$--invariant 
curves so that the image in $V\times S^1$ of each connected component 
is either a cylinder bounding a pair of 
type I Reeb orbits. 
The closure of the image of 
$u'$ in $V$ is a connected tree formed by trajectories 
of a gradient-like vector field. This tree contains a 
(perhaps broken) trajectory with endpoints 
$p_{\pm}=\pi (\gamma _{\pm})$. Moreover, $p_{\pm}$ are vertices 
of valent 1 of the tree, here the valent of a vertex is the 
number of edges coming out from this vertex as an endpoint. 
Suppose that the this trajectory 
contains other critical points. Then there is a critical point 
$p\neq p_{\pm}$ such that a trajectory between $p$ and $p_+$ is 
contained in the said broken trajectory from $p_-$ to $p_+$. 
The corresponding preimage of $p$ in $V\times S^1$ is the Reeb 
orbit $\gamma :=\gamma _p^m$. Since we must have the action 
inequalities 
\[ 
{\cal A}(\gamma _-)<{\cal A}(\gamma )<{\cal A}(\gamma _+), 
\] 
and since $\gamma$, $\gamma _-$, $\gamma _+$ have the same 
multiplicity we have 
\[ 
\text{ind}(p_-)>\text{ind}(p)>\text{ind}(p_+) .
\] 
But $\text{ind}(p_-)-\text{ind}(p_+)=1$, there exists no such $p$. 
So the trajectory between $p_{\pm}$ is unbroken and hence 
is equal to the tree. 

Now that the projection 
of the image of $u'$ in $V$ is a trajectory of a gradient-like 
vector field connecting critical points $p_-=\pi (\gamma _-)$ to  
$p_+=\pi (\gamma _+)$, so the linearized operator $D_{u'}$ is
surjective for generic $J'$, hence $u'$ is an isolated element of 
${\cal M}_{J'}(M';\gamma_-,\gamma_+)$. Thus $u'_{k_{\nu}}$ 
and hence $v'_{k_{\nu}}$ have 
to be $S^1$--invariant for all $\nu$ large enough, which contradicts 
with the assumption that there are non-$S^1$--invariant $u'_{k_{\nu}}$ 
for infinitely many $k_{\nu}$. Hence the lemma holds. 
\end{proof}

\begin{lem}
Let $(\alpha ',J')$ be a regular pair.
For $p\in Crit(f)$ with $\text{ind}(p)=k$, denote by
$\gamma ^m_p\in C_{2(n+m-1)-k}(\alpha ')$ the corresponding
principal generator with multiplicity $m\leq m_o$. Then
\[
\partial \gamma ^m_p = m\sum _{\text{ind}(q)=k+1}
\frac{a_q}{m}\gamma ^m_q,
\]
where $a_q$ is the algebraic number of trajectories of $\bar{J}X_h$ 
running from $q$ and $p$.
\end{lem}

Now let $m_o\to \infty$. After an easy computation on index  
we obtain the following: 

\begin{theo}   \label{HCM'}
Let $(M',\xi ')$ be a stabilization of a $(2n+1)$--dimensional
subcritical Stein-fillable contact manifold $(M,\xi )$, $n>2$,
and $(V,J)$ a subcritical Stein-filling of $(M,\xi )$. Then 
\[ 
HC_i(M',\xi')\cong \underset{m\in \mathbb{N}\cup \{ 0\} }
{\oplus}H_{2(n+m)-i}(V).  
\]  
\end{theo} 

Combining Theorem \ref{HCM'} with Theorem \ref{HCM=HCM'} we have the
following:

\begin{theo}  \label{HCM}
Let $(M,\xi )$ be a $(2n-1)$--dimensional subcritical Stein-fillable
contact manifold with $n>2$,
and $V$ a subcritical Stein-filling of $(M,\xi )$. Then 
\[ 
HC_i(M,\xi )\cong \underset{m\in \mathbb{N}\cup \{ 0\} }
{\oplus}H_{2(n+m-1)-i}(V).  
\] 
\end{theo}

When $n=2$ $(M,\xi )$ is the union of $S^3$ (3--dimensional contact 
0-handle) and a finite number of 3--dimensional contact 1--handles 
diffeomorphic to $\mathbb{R}\times S^2$. We write $s$ for the number 
of contact 1--handles of $M$. Let $\gamma _0$ denote the principal 
Reeb orbit in the 0-handle, and $\gamma _1$, $\gamma _2$, ..., 
$\gamma _s$ the principal Reeb orbits in each of the $s$' contact 
1--handles. These 1--handles can be attached to $S^3$ pairwise 
disjoint. Recall that when $n=2$ $c_1(\xi )=0$ so the 
$\bar{\mu}$--index of contractible Reeb orbits are independent of 
the spanning discs and hence are well-defined. 

Let $m_o\gg 1$ be a fixed positive integer, then 
 by deforming 0- and 1--handles we can obtain a 
suitable regular contact 1--form $\alpha$ and assume the 
following:  

\begin{cond} \label{cz-act} 
\begin{align} 
\overline{\gamma_j^m} & =\left\{ \begin{array}{ll} 
2m & \qquad j=0, \\ 2m-1 & \qquad j=1,2, ..., s  
\end{array}  \right. 
\quad \text{for } m\leq m_o  \label{m_o}   \\ 
{\cal A}(\gamma _1) & ={\cal A}(\gamma _2)= \cdots = 
{\cal A}(\gamma _s) \ll {\cal A}(\gamma _0)  \label{actions} 
\end{align}  
\end{cond} 

Let $J$ be a regular $\alpha$--admissible almost complex structure. 

\begin{prop}  \label{dimM} 
Assume ${\cal M}:={\cal M}_{J}(\gamma _-,\gamma _+)$ is not empty and 
$\overline{\gamma}_{\pm}\leq 2m_o$. 
\begin{enumerate} 
\item If $\dim {\cal M}=1$ then $\gamma _-=\gamma _j^m$, $\gamma _+
=\gamma _0^m$ for some $m\leq m_o$, $1\leq j\leq s$. 
\item If $\dim {\cal M}=2$ then $\gamma _-=\gamma _o^{m-1}$, $\gamma _+
=\gamma _0^m$ for some $m\leq m_o$. 
\end{enumerate} 
\end{prop}

Note that $C_1(\alpha )$ is nontrivial, it is generated by 
$\gamma _j$, $1\leq j\leq s$. Nevertheless Proposition \ref{dimM} 
and Condition \ref{cz-act} together imply that for $1<m\leq m_o$, 
the boundary of ${\cal M}(\gamma _0^{m-1},\gamma _0^m)$ does not 
contain any element of ${\cal M}(\gamma _j^m,\gamma _0^m)$. Also 
when $m=1$ the boundary of ${\cal M}(\gamma _0)$, where ${\cal M}
(\gamma _0)$ consists of holomorphic planes converging exponentially 
to $\gamma _0$ at $t=\infty$ at $z=\infty$, contains no 
holomorphic curves with more than one negative ends. 
Moreover,  we have the following result: 

\begin{lem} 
The boundary operator $\partial \co C_*(\alpha )\to C_{*-1}(\alpha )$ 
satisfies $\partial \circ \partial =0$, at least when $*\leq 2m_o$. 
\end{lem} 

\begin{proof} 
Condition \ref{cz-act} implies that for $1<m\leq m_o$, 
\[ 
\partial \gamma _j^m =0 \quad \forall  j=1,2, ...,s, 
\] 
and hence 
\[ 
\partial ^2\gamma _0^m=\partial (\sum _{j=1}^sa_j\gamma _j^m)=0 
\qquad \forall  j=1,2, ...,s.  
\] 
Now consider the $m=1$ case. Observe that $\partial ^2 (\gamma _j)=0$ 
for $j=1,...,s$ because $C_{-1}(\alpha )=0$. 
Also $\partial ^2 (\gamma _0)=0$ since 
$\partial \gamma _o =\sum_{j=1}^sb_j\gamma _j$ and 
$\partial \gamma _j=0$ for $j=1,...,s$. 
Thus $\partial ^2=0$ at least on $C_*(\alpha )$ with $*\leq 2m_o$. 
\end{proof} 

Therefore we can apply to $M$ the stabilization technique as before 
and obtain the $n=2$ version of Theorem \ref{HCM}. 

\begin{theo} 
Let $(M,\xi )$ be a 3--dimensional subcritical Stein-fillable
contact manifold, 
and $(V,J)$ a subcritical Stein domain such
that $\partial V=M$ and $\xi$ is the maximal complex
subbundle of $TM$. Then 
\[ 
HC_i(M,\xi)\cong \underset{m\in \mathbb{N}\cup \{ 0\} }
{\oplus}H_{2(n+m-1)-i}(V).  
\] 
\end{theo} 

This completes the proof of the Main Theorem.

\subsection*{Acknowledgements} 

The author is very grateful to Y Eliashberg for years of guidance, and
for many valuable discussion and suggestions.  The author is also
deeply thankful to S\,S Kim for reading the draft of this paper, and
to the referees for pointing out mistakes and giving many valuable
comments and suggestions on an earlier version of this paper.

\end{document}